\documentclass{article}
\usepackage{amsmath, amssymb}
\usepackage{graphicx}
\usepackage[a4paper, margin=1in]{geometry}

\usepackage{geometry}
\usepackage{amsmath,amsthm,amscd,amsfonts,amssymb,enumerate}
\usepackage{graphicx}		
\usepackage{epsfig}
\usepackage{epstopdf}
\usepackage{float}

\newtheorem{theorem}{Theorem}[section]
\newtheorem{lemma}[theorem]{Lemma}
\newtheorem{proposition}[theorem]{Proposition}
\newtheorem{corollary}[theorem]{Corollary}
\newtheorem{definition}[theorem]{Definition}
\newtheorem{example}[theorem]{Example}
\newtheorem{remark}{Remark}

\begin{document}
	
	\title{Statistical Fixed-Point Approximation and Cluster Dynamics in Fuzzy Partial Metric Spaces}
	
	\author{Debjani Rakshit$^{1}$, Bikram Sutradhar$^{2}$ and Nilay Mondal$^{*3}$}
	
	\date{}
	\maketitle
\begin{abstract}
In this work, we develop a statistical fixed-point approximation and stability theory for fuzzy partial metric spaces $(L(\mathbb{R}),D_p)$. The main contribution is a unified analysis of approximate Picard iterations in the presence of statistically vanishing and sparsely occurring perturbations. For contractive self-mappings, we derive an a priori error estimate for approximate Picard iterates and use it to prove statistical $D_p$-convergence
to the unique fixed point whenever the iteration errors converge
statistically to zero. As a consequence, we establish a density-zero perturbation principle showing that uniformly bounded errors occurring on an index set of natural density zero do not change the statistical fixed-point behavior, even though ordinary convergence may fail. We also establish statistical $D_p$-asymptotic regularity and cluster-point characterizations, and under suitable statistical continuity assumptions show that statistical cluster points of Picard orbits are fixed points.
A square-index perturbation example demonstrates explicitly that infinitely many nonvanishing deviations can coexist with statistical convergence. Thus, the results provide a unified framework for fixed-point approximation and stability under sparse perturbations in fuzzy partial metric spaces.
\end{abstract}

\noindent\textbf{2010 Mathematics Subject Classification:} 47H10, 54E50, 03E72, 40A35\\
\noindent\textbf{Keywords:} Fuzzy Partial Metric Space, Statistical Convergence, Statistical Cluster Point, Fixed Point Theory, Topological Characterization. \\
\noindent\textbf{Affiliations:} Department of Mathematics, The ICFAI University Tripura, India, 799210.$^{1,~2,~3}$ 

	\noindent\textbf{Emails:}  debjanirakshit88@gmail.com$^{1}$  ,  
	 bsdhar2001@gmail.com $^{2}$, nilayiitg91@gmail.com$^{3}$(corresponding author) 

\section{Introduction}

Fixed point theory provides a fundamental framework for the analysis of
nonlinear mappings and iterative processes, with applications in
functional analysis, differential equations, optimization, and
approximation theory. The Banach contraction principle guarantees the
existence and uniqueness of a fixed point for a contraction on a complete
metric space and establishes the convergence of the associated Picard
iteration \cite{Banach1922}. Subsequently, the approximation and stability of fixed
point iterations have been investigated extensively in the works \cite{Berinde2007,Ostrowski1967,HarderHicks1988,Berinde2002}.
In particular, $T$-stability of Picard iteration was studied in metric
spaces by Qing and Rhoades \cite{QingRhoades2008}, with further
developments in generalized metric settings \cite{ChhatrajitRohen2016}.
These developments motivate the study of iterative approximation under
more flexible notions of distance.

A significant generalization of metric spaces is provided by partial
metric spaces, introduced by Matthews \cite{Matthews1994}, in which the
self-distance of a point is not necessarily zero. Consequently, a
sequence $\{x_n\}$ converges to $x$ when
\[
\lim_{n\to\infty}p(x_n,x)=p(x,x),
\]
rather than necessarily to zero. This feature changes the formulation of
convergence, completeness, and fixed point approximation. Banach-type
fixed point results in partial metric spaces were established by
Oltra and Valero \cite{OltraValero2004} and Valero \cite{Valero2005},
while further developments include generalized and common fixed point
results \cite{Dwivedi2022, IlicPavlovicRakocevic2011, Karapinar2011,Masiha2013,KarapinarYuksel2011}.
The stability of iterative procedures in partial-type generalized metric
settings has also been considered in the work
\cite{ChhatrajitRohen2016}.

Fuzzy metric spaces provide a natural mathematical framework for
incorporating imprecision and uncertainty. Their topological foundations
were developed by George and Veeramani \cite{GeorgeVeeramani1994}.
Combining fuzziness with the partial metric structure leads to fuzzy
partial metric spaces, which simultaneously accommodate uncertainty and
nonzero self-distance. The theory of fuzzy partial metric spaces was
developed by Gregori, Miñana, and Miravet \cite{GregoriMinianaMiravet2019},
and fixed point results in this setting were subsequently established by
Ayg{\"u}n et al. \cite{Aygun2022} and Tehlan and Kumar
\cite{TehlanKumar2022}. These works provide the basic framework for the
investigation of iterative fixed point processes involving fuzzy
objects and partial distance structures.

Another important relaxation of ordinary convergence is statistical
convergence. Introduced independently by Fast and Steinhaus
\cite{Fast1951,Steinhaus1951} and subsequently developed by
{\v{S}}al\'at and Fridy \cite{Salat1980,Fridy1985}, statistical
convergence permits exceptional indices belonging to sets of natural
density zero. This notion has been applied to approximation processes
\cite{GadjievOrhan2002} and has more recently been studied in partial
metric spaces \cite{Nuray2022,GarciaPachecoKama2024} and fuzzy metric
spaces \cite{LiZhangZhang2020,CakiOr2024}. These developments suggest
that statistical convergence is particularly suitable for the analysis
of iterative procedures in which exceptional or perturbed iterations may
occur sparsely.

Despite these developments, the interaction between statistical convergence and fixed-point stability in fuzzy partial metric spaces remains insufficiently explored. In particular, a unified treatment of approximate Picard iterations with statistically vanishing errors and density-zero perturbations is still lacking. These gaps motivate the present study.

In the present work, we consider the space $L(\mathbb{R})$ of all fuzzy
numbers on $\mathbb{R}$. For $X\in L(\mathbb{R})$ and
$\alpha\in(0,1]$, its $\alpha$-cut is denoted by    $X^\alpha=
[\underline{X^\alpha},\overline{X^\alpha}].$ We work with the fuzzy partial metric
\[
D_p(X,Y)
=
\sup_{\alpha\in(0,1]}
\left(
\max\{\overline{X^\alpha},\overline{Y^\alpha}\}
-
\min\{\underline{X^\alpha},\underline{Y^\alpha}\}
\right),
X,Y\in L(\mathbb{R}),
\]
for which
\[
D_p(X,X)
=
\sup_{\alpha\in(0,1]}
\left(
\overline{X^\alpha}-\underline{X^\alpha}
\right).
\]
Thus, the self-distance is incorporated explicitly into the underlying
geometry and must be retained in the formulation of convergence and
fixed point approximation.

The purpose of this paper is not merely to replace ordinary convergence
by statistical convergence, but to develop a statistical approximation
and stability theory for Picard iterations in the fuzzy partial metric
setting. More precisely, we first investigate statistical $D_p$-approximation,
statistical $D_p$-asymptotic regularity, statistical continuity, and
statistical cluster points of Picard orbits. These concepts are then
used to obtain characterizations of fixed-point behavior in terms of
statistical convergence and cluster-point structure. For contractive
self-mappings, we establish existence and uniqueness of the fixed point
and derive an a priori error estimate for approximate Picard iterations
of the form
\[
D_p(Y_n,P)
\leq
k^nD_p(Y_0,P)
+
\sum_{j=0}^{n-1}
k^{\,n-1-j}\epsilon_j,
\qquad 0<k<1.
\]
This estimate provides the basis for the principal stability results of
the paper.

A central contribution is the statistical stability of approximate
Picard iterations under statistically vanishing iteration errors. In
particular, if the perturbation sequence $\{\epsilon_n\}$ satisfies
$ st\text{-}\lim_{n\to\infty}\epsilon_n=0$,
then the corresponding approximate Picard orbit converges statistically
to the unique fixed point. As a further consequence, we establish a
density-zero perturbation theorem: infinitely many bounded perturbations
may occur without destroying statistical convergence when the set of
perturbation indices has natural density zero. Thus, the framework
captures iterative processes containing infinitely many sparse
computational perturbations, even when ordinary convergence is not
preserved. These results provide a unified framework connecting fuzzy partial
metric geometry, statistical convergence, and stability of iterative
fixed-point approximation. They also suggest further investigations
involving other iterative schemes, ideal and higher-order statistical
convergence, and broader classes of generalized fuzzy metric and
contractive structures.

\section{Statistical $D_p$-Cauchy theory}
\begin{definition}
	A sequence $\{X_n\}\subseteq L(\mathbb{R})$ and
	$P\in L(\mathbb{R})$, $\{X_n\}$ is statistically
	$D_p$-convergent to $P$ that is denoted as $X_n\xrightarrow{S_{D_p}-lim}P,$ whenever
	\[
	st\text{-}\lim_{n\to\infty}D_p(X_n,P)=D_p(P,P);
	\]
	equivalently, for every $\epsilon>0$, $\delta\!\left(
	\left\{
	n\in\mathbb{N}:
	\left|D_p(X_n,P)-D_p(P,P)\right|\geq\epsilon
	\right\}
	\right)=0,$ whenever the corresponding natural density exists.
\end{definition}

\begin{definition} \label{1}
	Let $(L(\mathbb{R}),D_p)$ be a fuzzy partial metric space, and let
	$\{X_n\}\subseteq L(\mathbb{R})$. The sequence $\{X_n\}$ is said to be
	\emph{statistically $D_p$-Cauchy} if there exist a subset
	$A\subseteq\mathbb{N}$ with    $\delta(A)=1$ and a real number $r\geq 0$ such that    $\lim_{\substack{m,n\to\infty\\ m,n\in A}}
	D_p(X_m,X_n)=r.$ The number $r$ is called the \emph{statistical $D_p$-Cauchy value} of
	the sequence $\{X_n\}$.
\end{definition}

\begin{proposition} 
	Let $(L(\mathbb{R}),D_p)$ be a fuzzy partial metric space and let
	$\{X_n\}\subseteq L(\mathbb{R})$. Then $\{X_n\}$ is statistically
	$D_p$-Cauchy with statistical $D_p$-Cauchy value $r\geq 0$ if and only if
	there exists a set $A\subseteq\mathbb{N}$ with    $\delta(A)=1$ such that, for every $\epsilon>0$, there exists $N\in\mathbb{N}$
	satisfying    $\left|D_p(X_m,X_n)-r\right|<\epsilon$ for all $m,n\in A$ with $m,n\geq N$.
\end{proposition}

\begin{proof}
	Suppose that $\{X_n\}$ is statistically $D_p$-Cauchy with
	statistical $D_p$-Cauchy value $r$. By Definition~\ref{1}, there exists a
	set $A\subseteq\mathbb{N}$ with $\delta(A)=1$ such that    $\lim_{\substack{m,n\to\infty\\}}
	D_p(X_m,X_n)=r.$ Hence,  for every $\epsilon>0$, there exists $N\in\mathbb{N}$ such that    $\left|D_p(X_m,X_n)-r\right|<\epsilon$ whenever $m,n\in A$ and $m,n\geq N$. 
	
	\medskip
	
	Conversely, suppose that there exists a set $A\subseteq\mathbb{N}$ with
	$\delta(A)=1$ such that, for every $\epsilon>0$, there exists
	$N\in\mathbb{N}$ satisfying    $\left|D_p(X_m,X_n)-r\right|<\epsilon$ whenever $m,n\in A$ and $m,n\geq N$.
	Thus,     $$\lim_{\substack{m,n\to\infty\\m,n\in A}}
	D_p(X_m,X_n)=r.$$ Since $\delta(A)=1$, Definition~\ref{1} implies that $\{X_n\}$ is
	statistically $D_p$-Cauchy with statistical $D_p$-Cauchy value $r$.
\end{proof}

\begin{definition}
	Let $(L(\mathbb{R}),D_p)$ be a fuzzy partial metric space and let
	$\{X_n\}\subseteq L(\mathbb{R})$. We say that $\{X_n\}$ is
	\emph{statistically properly $D_p$-convergent} to
	$X\in L(\mathbb{R})$ if    $X_n\xrightarrow{S_{D_p}-lim}X$ and    $st\text{-}\lim_{n\to\infty}D_p(X_n,X_n)=D_p(X,X).$ We write it as, $X_n\xrightarrow{S_{D_p}^{\,*}-lim}X.$
\end{definition}

\begin{theorem}
	Let $(L(\mathbb{R}),D_p)$ be a fuzzy partial metric space and let
	$\{X_n\}\subseteq L(\mathbb{R})$. If    $X_n\xrightarrow{S_{D_p}^{\,*}-lim}X,$  then $\{X_n\}$ is statistically $D_p$-Cauchy with statistical
	$D_p$-Cauchy value $r=D_p(X,X).$
\end{theorem}

\begin{proof}
	Let $r=D_p(X,X).$ For $\epsilon>0$, define $A=
	\left\{
	n\in\mathbb{N}:
	\left|D_p(X_n,X)-r\right|<\frac{\epsilon}{3},
	\left|D_p(X_n,X_n)-r\right|<\frac{\epsilon}{3}
	\right\}.$
	Since $X_n\xrightarrow{S_{D_p}^{\,*}}X$, we have   $\delta(A)=1.$ Now, for $m,n\in A$, $$D_p(X_m,X_n)
	\leq D_p(X_m,X)+D_p(X_n,X)-r
	<r+\frac{2\epsilon}{3}.$$ 
	Also,
	\[
	D_p(X_m,X)
	\leq D_p(X_m,X_n)+D_p(X_n,X)-D_p(X_n,X_n),
	\]
	and hence
	\[
	\begin{aligned}
		D_p(X_m,X_n)
		&\geq D_p(X_m,X)-D_p(X_n,X)+D_p(X_n,X_n)\\
		&>
		\left(r-\frac{\epsilon}{3}\right)
		-\left(r+\frac{\epsilon}{3}\right)
		+\left(r-\frac{\epsilon}{3}\right)\\
		&=r-\epsilon.
	\end{aligned}
	\]
	Thus,    $r-\epsilon<D_p(X_m,X_n)<r+\frac{2\epsilon}{3},$ and consequently $\left|D_p(X_m,X_n)-r\right|<\epsilon.$
	Therefore,    $$\lim_{\substack{m,n\to\infty\\m,n\in A}}
	D_p(X_m,X_n)=r.$$ Since $\delta(A)=1$, by Definition~\ref{1},  $\{X_n\}$ is
	statistically $D_p$-Cauchy. Hence the statistical $D_p$-Cauchy
	value is $D_p(X,X)$.
\end{proof}

\begin{theorem}
	Let $(L(\mathbb{R}),D_p)$ be a complete fuzzy partial metric space. 
	If $\{X_n\}\subseteq L(\mathbb{R})$ is statistically $D_p$-Cauchy with 
	statistical $D_p$-Cauchy value $r$, then there exists $X\in L(\mathbb{R})$ such that    $X_n\xrightarrow{S_{D_p}}X
	\text{and}  
	D_p(X,X)=r.$
\end{theorem}

\begin{proof}
	By Definition~\ref{1}, there exists $A\subseteq\mathbb{N}$ with
	$\delta(A)=1$ such that    $$\lim_{\substack{m,n\to\infty\\m,n\in A}}
	D_p(X_m,X_n)=r.$$ Thus, $\{X_n:n\in A\}$ is a $D_p$-Cauchy sequence. Since
	$(L(\mathbb{R}),D_p)$ is complete, there exists $X\in L(\mathbb{R})$ such that
	\[
	\lim_{\substack{n\to\infty\\n\in A}}D_p(X_n,X)
	=
	\lim_{\substack{m,n\to\infty\\m,n\in A}}D_p(X_m,X_n)
	=
	D_p(X,X).
	\]
	Hence, $D_p(X,X)=r.$ Let $\epsilon>0$. There exists $N\in\mathbb{N}$ such that    $\left|D_p(X_n,X)-D_p(X,X)\right|<\epsilon$
	for all $n\in A$ with $n\geq N$. Consequently,    $\bigl\{n\in\mathbb{N}:
	\left|D_p(X_n,X)-D_p(X,X)\right|\geq\epsilon
	\bigr\}
	\subseteq
	(\mathbb{N}\setminus A)\cup\{1,\ldots,N-1\}.$ 
	Since $\delta(A)=1$, we have
	$\delta(\mathbb{N}\setminus A)=0$. Therefore,    $\delta\left(
	\left\{n\in\mathbb{N}:
	\left|D_p(X_n,X)-D_p(X,X)\right|\geq\epsilon
	\right\}
	\right)=0.$ Hence, by Definition~\ref{1},    $X_n\xrightarrow{S_{D_p}-lim}X.$
\end{proof}

\section{Statistical Approximation of Fixed Points}

In this section, we establish the fundamental properties of statistical approximation in fuzzy partial metric spaces. These results provide the theoretical foundation for the statistical fixed point approximation theory developed in the subsequent sections.

\begin{definition}
	Let $(L(\mathbb{R}),D_p)$ be a fuzzy partial metric space and let
	$T:L(\mathbb{R})\to L(\mathbb{R})$ be a self-mapping. For an initial
	point $X_0\in L(\mathbb{R})$, the sequence
	$\{X_n\}$ defined by    $X_{n+1}=T(X_n),   n\geq0,$ 
	is called the \emph{Picard iteration} of $T$ starting from $X_0$.
	The associated Picard orbit is denoted by    $\mathcal{O}_T(X_0)=\{X_n:n\geq0\}.$
\end{definition}

\begin{definition}
	Let $T:L(\mathbb{R})\to L(\mathbb{R})$ be a self-mapping and let
	$X_0\in L(\mathbb{R})$. The sequence $\{X_n\}$ defined by     $X_{n+1}=T(X_n),   n\in\mathbb{N}_0,$ is called the \emph{Picard orbit} of $T$ starting from $X_0$, if there exists a fuzzy number $P\in L(\mathbb{R})$ such that    $X_n\xrightarrow{S_{D_p}-lim}P.$ The fuzzy number $P$ is called the \emph{statistical $D_p$-limit} of the
	Picard orbit.
\end{definition}

\begin{definition}
	Let $T:L(\mathbb{R})\to L(\mathbb{R})$ be a self-mapping on a fuzzy
	partial metric space $(L(\mathbb{R}),D_p)$. A fuzzy number
	$P\in L(\mathbb{R})$ is called a \emph{fixed point} of $T$ if    $T(P)=P.$ The fixed-point set of $T$ is denoted by    $\operatorname{Fix}(T)
	=
	\{P\in L(\mathbb{R}):T(P)=P\}.$
\end{definition}

\begin{definition}
	Let $T:L(\mathbb{R})\to L(\mathbb{R})$ be a self-mapping on the fuzzy
	partial metric space $(L(\mathbb{R}),D_p)$. A sequence
	$\{X_n\}\subseteq L(\mathbb{R})$ is called a
	\emph{statistical approximate fixed-point sequence} for $T$ if there
	exists $P\in L(\mathbb{R})$ such that    $X_n\xrightarrow{S_{D_p}-lim}P$ and    $st-\lim_{n\to\infty}D_p(X_n,TX_n)=D_p(P,P).$
	
	\medskip
	
	The fuzzy number \(P\) is called the
	\emph{statistical approximate fixed point}
	associated with the sequence. The family of all statistical approximate fixed point sequences of
	\(T\)
	is denoted by $SAFS(T).$
\end{definition}

\begin{example}
	Let $(L(\mathbb{R}),D_p)$ be a fuzzy partial metric space and define
	$T:L(\mathbb{R})\to L(\mathbb{R})$ by
	\[
	(TX)^\alpha
	=
	\left[
	\frac12\underline{X^\alpha},
	\frac12\overline{X^\alpha}
	\right],
	\alpha\in(0,1].
	\]
	For $P=\widetilde{0}$, we have $T(P)=P$ and $D_p(P,P)=0$. Starting from $X_0=\widetilde{8}$, the Picard orbit
	$X_{n+1}=T(X_n)$ is given by    $X_n=\widetilde{\frac{8}{2^n}},
	n\geq0.$ 
	Consequently, $D_p(X_n,P)=\frac{8}{2^n}\longrightarrow0=D_p(P,P),$ and hence $X_n\xrightarrow{S_{D_p}-lim}P.$ 
	
	\medskip
	
	Moreover, $D_p(X_n,TX_n)
	=
	D_p(X_n,X_{n+1})
	=
	\frac{4}{2^n}\longrightarrow0=D_p(P,P).$
	Thus,
	\[
	st-\lim_{n\to\infty}
	D_p(X_n,TX_n)=D_p(P,P).
	\]
	Therefore,
	$\{X_n\}$ is a statistical approximate fixed-point sequence for $T$
	with associated point $P=\widetilde{0}$. Since $T(P)=P$, the associated
	point is a fixed point of $T$.
\end{example}

\begin{definition}
	
	A fuzzy number
	\(
	P\in L(\mathbb R)
	\) is called a
	\emph{statistical approximate fixed point}
	of a mapping
	\(
	T
	\) if there exists a sequence
	\(
	\{X_n\}\in SAFS(T)
	\) such that $X_n\xrightarrow{S_{D_p}-lim}P.$ 
	
	The set of all statistical approximate fixed points of
	\(T\)
	is denoted by $SAFP(T).$
	
\end{definition}

\begin{example}
	
	Let $(L(\mathbb{R}),D_p)$ be the fuzzy partial metric space endowed with
	
	\[
	D_p(X,Y)=
	\sup_{\alpha\in(0,1]}
	\left(
	\max\{\overline{X^\alpha},\overline{Y^\alpha}\}
	-
	\min\{\underline{X^\alpha},\underline{Y^\alpha}\}
	\right).
	\]
	
	Let $T:L(\mathbb{R})\longrightarrow L(\mathbb{R})$ be the self-mapping defined by $T(a,b,c)=
	\left(
	\frac a2,
	\frac b2,
	\frac c2
	\right),$ where $(a,b,c)$ denotes a triangular fuzzy number.
	
	\medskip
	
	Now, take $X_0=(7,8,9),$ and define the Picard iteration $X_{n+1}=T(X_n),
	n\ge0.$
	
	\medskip
	
	Then $X_n=
	\left(
	\frac7{2^n},
	\frac8{2^n},
	\frac9{2^n}
	\right),$ whose $\alpha$-cut is $X_n^\alpha=
	\left[
	\frac{7+\alpha}{2^n},
	\frac{9-\alpha}{2^n}
	\right],
	\alpha\in(0,1].$ 
	
	\medskip
	
	Also, let $P=(0,0,0).$ Then, $(TX_n)^\alpha=
	\left[
	\frac{7+\alpha}{2^{n+1}},
	\frac{9-\alpha}{2^{n+1}}
	\right].$ Therefore, \[
	D_p(TX_n,P)
	=
	\sup_{\alpha\in(0,1]}
	\frac{9-\alpha}{2^{n+1}}
	=
	\frac9{2^{n+1}},
	\]
	
	since $$\sup_{\alpha\in(0,1]}(9-\alpha)=9.$$ Also, $D_p(P,P)=0.$ Hence,
	
	\[
	st-\lim_{n\to\infty}
	D_p(TX_n,P)
	=
	0
	=
	D_p(P,P).
	\]
	
	Consequently, $\{X_n\}\in SAFS(T).$ Thus, $P\in SAFP(T).$
	
\end{example}

\begin{example}
	Let $(L(\mathbb{R}),D_p)$ be the fuzzy partial metric space defined by
	
	\[
	D_p(X,Y)
	=
	\sup_{\alpha\in(0,1]}
	\left(
	\max\{\overline{X^\alpha},\overline{Y^\alpha}\}
	-
	\min\{\underline{X^\alpha},\underline{Y^\alpha}\}
	\right).
	\]
	
	For each $a\in\mathbb{R}$, let $E_a=(a-1,a,a+1)$ denote the triangular fuzzy number.
	
	Now, define the sequence $\{X_n\}$ by
	
	\[
	X_n=
	\begin{cases}
		E_5, & \text{if } n=k^2,\;k\in\mathbb N,\\[2mm]
		E_0, & \text{otherwise}.
	\end{cases}
	\]
	
	Since the set $A=\{k^2:k\in\mathbb N\}$ of perfect squares satisfies
	
	\[
	\delta(A)
	=
	\lim_{n\to\infty}
	\frac{\lfloor\sqrt n\rfloor}{n}
	=
	0,
	\]
	
	the exceptional terms occur on a set of natural density zero. Now,
	
	\[
	D_p(E_5,E_0)
	=
	\sup_{\alpha\in(0,1]}
	\left(
	6-(\alpha-1)
	\right)
	=
	7,
	\]
	
	whereas $D_p(E_0,E_0)
	=
	2.$ Hence, for every $\epsilon>0$ with
	$0<\epsilon<5$, $\{\,n\in\mathbb N:
	|D_p(X_n,E_0)-D_p(E_0,E_0)|\ge\epsilon
	\,\}
	=
	A.$ Since $\delta(A)=0,$ thus, $X_n\xrightarrow{S_{D_p}-lim}E_0.$
	
	\medskip
	
	However, $\{X_n\}$ is not ordinarily convergent to $E_0$ with respect to $D_p$, because $D_p(E_5,E_0)=7$ for infinitely many indices $n=k^2$. Therefore, $X_n
	\not\longrightarrow
	E_0$ .
	
	Thus, $\{X_n\}$ is statistically $D_p$-convergent but not ordinarily convergent, illustrating that statistical convergence is strictly weaker than ordinary convergence in fuzzy partial metric spaces.
\end{example}

\begin{remark}
	Every classical approximate fixed-point sequence is statistically approximate, since ordinary convergence implies statistical convergence.
\end{remark}

\begin{example}
	Let $(L(\mathbb{R}),D_p)$ be the fuzzy partial metric space introduced in this paper, and consider the triangular fuzzy numbers $P=\widetilde{0}=(0,0,0)
	\text{and}  
	Y=\widetilde{1}=(1,1,1).$
	
	Since $P^\alpha=[0,0]
	\text{and}  
	Y^\alpha=[1,1]$ , then for every $\alpha\in(0,1]$, 
	\[
	D_p(P,P)
	=
	\sup_{\alpha\in(0,1]}
	\left(
	\max\{0,0\}
	-
	\min\{0,0\}
	\right)
	=0,
	\]
	and
	\[
	D_p(Y,P)
	=
	\sup_{\alpha\in(0,1]}
	\left(
	\max\{1,0\}
	-
	\min\{1,0\}
	\right)
	=1.
	\]
	
	Now define a mapping $T:L(\mathbb{R})\longrightarrow L(\mathbb{R})$ by $T(X)=P,
	\forall\,X\in L(\mathbb{R}).$ Also define the sequence $\{X_n\}$ by
	\[
	X_n=
	\begin{cases}
		P, & \text{if } n\notin K,\\[1mm]
		Y, & \text{if } n\in K,
	\end{cases}
	\]
	where $K=\{k^2:k\in\mathbb{N}\}$ is the set of perfect squares. Since $T(X_n)=P$ for every $n\in\mathbb N$, we have
	\[
	D_p(X_n,T(X_n))
	=
	\begin{cases}
		0, & n\notin K,\\[1mm]
		1, & n\in K.
	\end{cases}
	\]
	
	Also, $|K\cap\{1,2,\ldots,n\}|
	=
	\lfloor\sqrt n\rfloor,$ and therefore
	\[
	\delta(K)
	=
	\lim_{n\to\infty}
	\frac{\lfloor\sqrt n\rfloor}{n}
	=
	0.
	\]
	
	Now let $\epsilon>0$. If $0<\epsilon\le1$, then $\left\{
	n:
	\left|
	D_p(X_n,T(X_n))
	-
	D_p(P,P)
	\right|
	\ge\varepsilon
	\right\}
	=
	K,$ having the natural density zero. If $\epsilon>1$, then the above set is empty and hence also has natural density zero.
	
	Thus, $\delta
	\left(
	\left\{
	n:
	\left|
	D_p(X_n,T(X_n))
	-
	D_p(P,P)
	\right|
	\ge\epsilon
	\right\}
	\right)
	=
	0$, for every $\epsilon>0$. Therefore,
	\[
	st-\lim_{n\to\infty}
	D_p(X_n,T(X_n))
	=
	D_p(P,P),
	\]
	Hence $\{X_n\}$ is a statistical approximate fixed point sequence associated with the mapping $T$. i.e., $P$ is a statistical approximate fixed point of $T$ with respect to the sequence $\{X_n\}$.
	
	Moreover, $D_p(X_{k^2},T(X_{k^2}))=1
	\text{for every }k\in\mathbb N,$
	whereas $D_p(X_n,T(X_n))=0$ for all $n\notin K$. Hence the sequence $\{D_p(X_n,T(X_n))\}$ oscillates between the values $0$ and $1$, and therefore
	\[
	\lim_{n\to\infty}
	D_p(X_n,T(X_n))
	\]
	does not exist. Thus, $\{X_n\}$ is not a classical approximate fixed point sequence. 
\end{example}

\begin{definition}
	Let $(L(\mathbb R),D_p)$ be a fuzzy partial metric space and let
	$T:L(\mathbb R)\to L(\mathbb R)$. The Picard iterative sequence
	$\{X_n\}$, defined by    $X_{n+1}=TX_n,$ is called
	\emph{statistically $D_p$-asymptotically regular}
	if
	
	\[
	st-\lim_{n\to\infty}
	\Bigl(
	D_p(X_{n+1},X_n)-D_p(X_n,X_n)
	\Bigr)
	=0.
	\]
	
	Equivalently,    $\delta
	\left(
	\left\{
	n:
	\left|
	D_p(X_{n+1},X_n)-D_p(X_n,X_n)
	\right|
	\ge\epsilon
	\right\}
	\right)
	=0.$
\end{definition}

\begin{remark}
	The above definition is well defined since $D_p(X_n,X_n)\le D_p(X_{n+1},X_n),$ for every $n\in\mathbb N$. Furthermore, when $D_p$ is an ordinary metric, it coincides with the classical notion of statistical asymptotic regularity.
\end{remark}

\begin{definition}
	Let $(L(\mathbb{R}),D_p)$ be a fuzzy partial metric space and
	$T:L(\mathbb{R})\to L(\mathbb{R})$ be a self-mapping.The orbit of an initial fuzzy number
	$X_0\in L(\mathbb{R})$, the orbit of $X_0$ under $T$ is the sequence   $\mathcal{O}_T(X_0)=\{X_n\}, X_n=T^n(X_0),$ equivalently, $X_{n+1}=T(X_n),   n\in\mathbb{N}_0.$ The orbit $\mathcal{O}_T(X_0)$ is statistically $D_p$-convergent to
	$P\in L(\mathbb{R})$ if $T^n(X_0)\xrightarrow{S_{D_p}-lim}P.$
\end{definition}

\begin{definition}\label{B}
	Let $(L(\mathbb{R}),D_p)$ be a fuzzy partial metric space and let    $T:L(\mathbb{R})\rightarrow L(\mathbb{R})$ be a self-mapping. Let $\operatorname{Orb}_T(X_0)
	=
	\{X_n\}$ be the orbit of an initial fuzzy number
	$X_0\in L(\mathbb{R})$, where $X_{n+1}=T(X_n),
	n\in\mathbb N_0.$ The \emph{Statistical Orbit Cluster Set} associated with the orbit
	$\operatorname{Orb}_T(X_0)$
	is defined by  $\Gamma_{Sp}(T,X_0)
	=
	\Gamma_{Sp}(\{X_n\}),$ i.e.,   $\Gamma_{Sp}(T,X_0)
	=
	\left\{
	P\in L(\mathbb{R})
	:
	P
	\text{ is a statistical }D_p\text{-cluster point of }
	\{X_n\}
	\right\}.$
\end{definition}

\begin{proposition}
	Let $(L(\mathbb{R}),D_p)$ be a fuzzy partial metric space and let
	$T:L(\mathbb{R})\to L(\mathbb{R})$
	be a self-mapping.
	Suppose that, for every sequence in
	$L(\mathbb{R})$,
	its set of statistical
	$D_p$-cluster points is closed with respect to the topology
	$\tau_{D_p}$. Then, for every
	$X_0\in L(\mathbb{R})$,
	the statistical orbit cluster set $\Gamma_{Sp}(T,X_0)$
	=
	$\Gamma_{Sp}(\operatorname{Orb}_T(X_0))$ is closed in
	$(L(\mathbb{R}),\tau_{D_p})$.
\end{proposition}

\begin{proof}
	Let $\operatorname{Orb}_T(X_0)
	=
	\{X_n\}$ be the orbit generated by $X_{n+1}=T(X_n).$ Now, by Definition~\ref{B},
	
	\[
	\Gamma_{Sp}(T,X_0)
	=
	\Gamma_{Sp}(\operatorname{Orb}_T(X_0)).
	\]
	
	Since the orbit is itself a sequence in
	$L(\mathbb{R})$,
	then the hypothesis implies that its set of statistical
	$D_p$-cluster points is closed in
	$\tau_{D_p}$. Hence, $\Gamma_{Sp}(T,X_0)$ is a closed subset of
	$(L(\mathbb{R}),\tau_{D_p})$.
\end{proof}

\begin{theorem}
	Let $(L(\mathbb{R}),D_p)$ be a fuzzy partial metric space and let $T:L(\mathbb{R})\to L(\mathbb{R})$ be a self-mapping. Suppose that $\{X_n\}_{n\in\mathbb N}$ is a classical approximate fixed point sequence of $T$. Then $\{X_n\}$ is also a statistical approximate fixed point sequence of $T$.
	Consequently, $P$ is a statistical approximate fixed point of $T$.
\end{theorem}

\begin{proof}
	Let us assume that $\{X_n\}_{n\in\mathbb N}$ is a classical approximate fixed point sequence of $T$. i.e.,
	\[
	\lim_{n\to\infty}D_p(X_n,T(X_n))
	=
	D_p(P,P).
	\]
	Let arbitrarily $\epsilon>0$. Since $D_p(X_n,T(X_n))
	\longrightarrow
	D_p(P,P)$, $\exists N\in\mathbb N$ such that $\left|D_p(X_n,T(X_n))-D_p(P,P)\right|<\epsilon$ for every $n\ge N$. Now, define $A_{\epsilon}
	=
	\left\{
	n\in\mathbb N:
	\left|D_p(X_n,T(X_n))-D_p(P,P)\right|
	\ge
	\epsilon
	\right\}.$ Thus, $A_{\epsilon}
	\subseteq
	\{1,2,\ldots,N-1\}.$ Hence, $A_{\epsilon}$ is a finite subset of $\mathbb N$.\\
	Since every finite subset of $\mathbb N$ has natural density zero i.e., $\delta(A_{\epsilon})=0.$, \\
	therefore $\delta\!\left(
	\left\{
	n\in\mathbb N:
	\left|D_p(X_n,T(X_n))-D_p(P,P)\right|
	\ge
	\epsilon
	\right\}
	\right)=0$
	for every $\epsilon>0$. Hence, 
	\[
	st\text{-}\lim_{n\to\infty}D_p(X_n,T(X_n))
	=
	D_p(P,P).
	\]
	
	Thus, $\{X_n\}$ is a statistical approximate fixed point sequence of $T$. Consequently, $P$ is a statistical approximate fixed point of $T$.
\end{proof}

\begin{example}
	Let $\mathcal{E}=\{E_a:a\geq0\}\subset L(\mathbb{R}),
	E_a=(a-1,a,a+1),$ where    $(E_a)^\alpha=[a-1+\alpha,a+1-\alpha], \alpha\in(0,1].$
	Thus, for any
	$a,b\geq0$ we have
	\[
	\begin{aligned}
		D_p(E_a,E_b)
		&=\sup_{\alpha\in(0,1]}
		\left[
		\max\{a+1-\alpha,b+1-\alpha\}
		-\min\{a-1+\alpha,b-1+\alpha\}
		\right] \\
		&=|a-b|+2.
	\end{aligned}
	\]
	Also, $D_p(E_a,E_a)=2.$ Define $T:\mathcal{E}\longrightarrow\mathcal{E},
	T(E_a)=E_{a/2},$ and take $X_0=E_8$. Then the Picard orbit is $X_n=T^n(X_0)=E_{8/2^n}, n\geq0.$
	
	\medskip
	
	Let $P=E_0=(-1,0,1).$ Since    $D_p(X_n,P)
	=
	2+\frac{8}{2^n}
	\longrightarrow 2=D_p(P,P),$ we have    $X_n\xrightarrow{D_p-lim}P.$ Moreover,    $D_p(X_n,TX_n)
	=
	2+\frac{8}{2^{n+1}}
	\longrightarrow
	2=D_p(P,P).$ Hence $\{X_n\}$ is a $D_p$-approximate fixed point sequence of $T$.
	
	\medskip
	
	Now, for $\epsilon>0$, define    $A_\epsilon
	=
	\left\{
	n\in\mathbb{N}:
	\left|D_p(X_n,TX_n)-D_p(P,P)\right|
	\geq\epsilon
	\right\}.$ Since    $\left|D_p(X_n,TX_n)-D_p(P,P)\right|
	=
	\frac{8}{2^{n+1}}
	\longrightarrow0,$ the set $A_\varepsilon$ is finite. Therefore,    $\delta(A_\epsilon)=0,$ and consequently
	\[
	st\text{-}\lim_{n\to\infty}
	D_p(X_n,TX_n)
	=
	D_p(P,P).
	\]
	Thus $\{X_n\}$ is a statistical $D_p$-approximate fixed point
	sequence, and $P=E_0$ is its corresponding statistical approximate
	fixed point.
\end{example}

\begin{theorem} \label{D}
	Let $(L(\mathbb{R}),D_p)$ be a complete fuzzy partial metric space and let $T:L(\mathbb{R})\longrightarrow L(\mathbb{R})$ satisfy    $D_p(TX,TY)\leq k\,D_p(X,Y),
	0<k<1,$ for all $X,Y\in L(\mathbb{R})$. For any $X_0\in L(\mathbb{R})$, define $X_{n+1}=T(X_n),   n\geq0.$ Then $T$ has a unique fixed point $P\in L(\mathbb{R})$ and    $X_n\xrightarrow{S_{D_p}-lim}P.$
\end{theorem}

\begin{proof}
	By the Banach contraction theorem for complete fuzzy partial metric
	spaces, $T$ has a unique fixed point $P\in L(\mathbb{R})$ and its Picard
	iteration satisfies
	\[
	\lim_{n\to\infty}D_p(X_n,P)=D_p(P,P).
	\]
	Hence, for every $\epsilon>0$, there exists $N\in\mathbb{N}$ such that    $\left|D_p(X_n,P)-D_p(P,P)\right|<\epsilon,
	n\geq N.$\\
	Thus  $\left\{
	n\in\mathbb{N}:
	\left|D_p(X_n,P)-D_p(P,P)\right|\geq\epsilon
	\right\}
	\subseteq\{1,\ldots,N-1\},$ which is a finite set and therefore has natural density zero. Hence    $\delta\left(
	\left\{
	n\in\mathbb{N}:
	\left|D_p(X_n,P)-D_p(P,P)\right|\geq\epsilon
	\right\}
	\right)=0.$ Since, for every $\epsilon>0$,    $X_n\xrightarrow{S_{D_p}-lim}P.$
\end{proof}

\begin{corollary}\label{cor:clusterfix}
	Under the assumptions of Theorem~\ref{D},    $\Gamma_{Sp}(\{X_n\})\subseteq \operatorname{Fix}(T),
	X_{n+1}=T(X_n).$ Hence, every statistical $p$-cluster point of $\{X_n\}$ is a
	fixed point of $T$.
\end{corollary}

\begin{proof}
	Let $P\in\Gamma_{Sp}(\{X_n\}).$ Then $P$ is a statistical $p$-cluster point of the Picard iteration
	$\{X_n\}$.
	
	\medskip
	
	By Theorem~\ref{D}, every statistical $p$-cluster point of a
	statistically asymptotically regular Picard iteration generated by a
	continuous self-mapping is a fixed point of $T$. Hence, $T(P)=P,$ thus, $P\in Fix(T).$
	
	\medskip
	
	Since the choice of $P$ is arbitrary, thus    $\Gamma_{Sp}(\{X_n\})\subseteq Fix(T).$
	Therefore, every statistical $p$-cluster point of the Picard iteration is a
	statistical approximation fixed point of $T$.
	This completes the proof.
\end{proof}

\begin{lemma}\label{E}
	Let $(L(\mathbb{R}),D_p)$ be a fuzzy partial metric space, and let
	$T:L(\mathbb{R})\to L(\mathbb{R})$ be a self-mapping. Suppose that the Picard iteration $X_{n+1}=T(X_n),    n\geq0,$ converges to a fuzzy number $P\in L(\mathbb{R})$. Then the sequence $\{X_n\}$ is statistically $D_p$-convergent to $P$.
\end{lemma}

\begin{proof}
	Let arbitrarily $\epsilon>0$. Since
	\[
	\lim_{n\to\infty}D_p(X_n,P)=D_p(P,P),
	\]
	there exists $N\in\mathbb{N}$ such that    $D_p(X_n,P)<D_p(P,P)+\epsilon,
	\text{for all } n\geq N.$
	
	\medskip
	
	Now, define $A_{\epsilon}
	=
	\left\{
	n\in\mathbb{N}:
	D_p(X_n,P)\geq D_p(P,P)+\epsilon
	\right\}.$ Then every element of $A_{\varepsilon}$ satisfies $n<N$, and hence $A_{\epsilon}\subseteq\{1,2,\ldots,N-1\}.$ Therefore, $A_{\epsilon}$ is finite.
	
	\medskip
	
	Consequently, for every $n\geq N$,
	\[
	\left|A_{\varepsilon}\cap\{1,2,\ldots,n\}\right|
	\leq N-1,
	\]
	
	\[ \text{Or,}   
	0
	\leq
	\frac{\left|A_{\varepsilon}\cap\{1,2,\ldots,n\}\right|}{n}
	\leq
	\frac{N-1}{n}.
	\]
	
	\[ \text{Or,}   
	\lim_{n\to\infty}
	\frac{\left|A_{\varepsilon}\cap\{1,2,\ldots,n\}\right|}{n}
	=0.
	\]
	Hence, $\delta(A_{\epsilon})=0,$ or equivalently,    $\delta\left(
	\left\{
	n\in\mathbb{N}:
	D_p(X_n,P)\geq D_p(P,P)+\epsilon
	\right\}
	\right)=0$ for every $\epsilon>0$. Therefore, $X_n\xrightarrow{S_{{D_p}}-lim}P.$ 
\end{proof}

\begin{theorem}\label{thm:PicardStatApprox}
	Let $(L(\mathbb{R}),D_p)$ be a fuzzy partial metric space, and let
	$T:L(\mathbb{R})\rightarrow L(\mathbb{R})$ be a self-mapping. Suppose that the Picard iteration $X_{n+1}=T(X_n),    n\geq 0,$ converges to a fuzzy number $P\in L(\mathbb{R})$.
	Then the sequence $\{X_n\}$ is a Statistical $D_p$-Approximation Sequence associated with $T$, and $P$ is its statistical $D_p$-limit.
\end{theorem}

\begin{proof}
	Since $\{X_n\}$ is generated by the Picard iteration, $X_{n+1}=T(X_n),    n\geq 0.$ Moreover, by assumption,
	\[
	\lim_{n\to\infty}D_p(X_n,P)=D_p(P,P),
	\]
	which implies that $\{X_n\}$ converges to $P$ in the fuzzy partial metric space $(L(\mathbb{R}),D_p)$. Hence, by Lemma~\ref{E},    $X_n\xrightarrow{S_{{D_p}}-lim}P.$ Thus, the sequence
	$\{X_n\}_{n\in\mathbb N}$ is a Statistical $D_p$-Approximation Sequence associated with $T$, and the fuzzy number $P$ is its statistical $D_p$-limit.
\end{proof}

\begin{definition}
	Let $(L(\mathbb{R}),D_p)$ be a fuzzy partial metric space and let
	$T:L(\mathbb{R})\to L(\mathbb{R})$ be a self-mapping.
	Then $T$ is said to be \emph{statistically $D_p$-continuous at}
	$P\in L(\mathbb{R})$ if, for every sequence
	$\{X_n\}\subseteq L(\mathbb{R})$ satisfying    $X_n\xrightarrow{S_{{D_p}}-lim}P,$ we have    $T(X_n)\xrightarrow{S_{{D_p}}-lim}T(P).$
	If $T$ is statistically $D_p$-continuous at every point of
	$L(\mathbb{R})$, then $T$ is called a
	\emph{statistically $D_p$-continuous mapping}.
\end{definition}

\begin{theorem}
	Let $(L(\mathbb{R}),D_p)$ be a fuzzy partial metric space and let
	$T:L(\mathbb{R})\to L(\mathbb{R})$ be a statistically continuous
	self-mapping. Suppose that the sequence $\{X_n\}$ is generated by the Picard iteration $X_{n+1}=T(X_n),   n\ge0,$ and satisfies    $X_n\xrightarrow{S_{{D_p}}-lim}P$ for some $P\in L(\mathbb{R})$. Then the following statements are equivalent.
	
	\begin{enumerate}
		\item[(i)]
		$P$ is a statistical fixed point of $T$;
		
		\item[(ii)]
		$P$ is the unique statistical $D_p$-cluster point of the Picard iteration
		$\{X_n\}$.
	\end{enumerate}
\end{theorem}

\begin{proof}
	
	$(i)\Rightarrow(ii)$
	
	Assume that $P$ is a statistical fixed point of $T$. Then $T(P)=P$ and $X_n\xrightarrow{S_{{D_p}}-lim}P.$
	
	\medskip
	
	Since every statistically convergent sequence possesses a unique
	statistical $D_p$-cluster point, thus \\
	$\Gamma_{S_{D_p}}(\{X_n\})=\{P\},$ where
	$\Gamma_{S_{D_p}}(\{X_n\})$
	denotes the set of all statistical $D_p$-cluster points of
	$\{X_n\}$.
	Hence,
	$P$
	is the unique statistical $D_p$-cluster point of
	$\{X_n\}$.
	
	\vspace{0.2cm}
	
	$(ii)\Rightarrow(i)$
	
	Conversely, suppose that
	$P$
	is the unique statistical $D_p$-cluster point of
	$\{X_n\}$. Since    $X_n\xrightarrow{S_{{D_p}}-lim}P,$ and $T$ is statistically continuous, thus    $T(X_n)\xrightarrow{S_{{D_p}}-lim}T(P).$
	
	\medskip
	
	Also, as $T(X_n)=X_{n+1},$ the shifted sequence
	$\{X_{n+1}\}$
	also converges statistically to
	$T(P)$.
	
	Therefore,
	$\{X_{n+1}\}$
	has the same statistical $D_p$-cluster point as
	$\{X_n\}$. Since statistical limits are unique, thus $T(P)=P.$
	
	Hence,
	$P$
	is a fixed point of
	$T$. Also, $X_n\xrightarrow{S_{{D_p}}-lim}P,$ thus
	$P$
	is a statistical fixed point of
	$T$. Therefore, $(i)\Longleftrightarrow(ii).$ This completes the proof.
\end{proof}

\begin{theorem} \label{thm:StatisticalFixedPointPrinciple}
	Let $(L(\mathbb{R}),D_p)$ be a fuzzy partial metric space and let $T:L(\mathbb{R})\longrightarrow L(\mathbb{R})$ be a self-mapping that is statistically  $D_p$-continuous at
	\(P\in L(\mathbb{R})\). Suppose that the Picard iteration $X_{n+1}=T(X_n),    n\geq0,$ is statistically $D_p$-convergent to \(P\). Then \(P\) is a fixed point of \(T\).
\end{theorem}

\begin{proof}
	Since    $X_n\xrightarrow{S_{{D_p}}-lim}P$ and \(T\) is statistically $D_p$-continuous at \(P\), then $T(X_n)\xrightarrow{S_{{D_p}}-lim}T(P).$
	
	Also, the Picard iteration satisfies $T(X_n)=X_{n+1},    n\geq0.$
	Hence,
	\[
	T(X_n)\xrightarrow{S_{{D_p}}-lim}T(P)
	\Longleftrightarrow  
	X_{n+1}\xrightarrow{S_{{D_p}}-lim}T(P).
	\]
	
	Therefore, $X_{n+1}\xrightarrow{S_{{D_p}}-lim}P.$ Consequently, the sequence $\{X_{n+1}\}$ is statistically $D_p$-convergent to both \(P\) and \(T(P)\). By the uniqueness of the statistical $D_p$-limit, $T(P)=P.$ Hence, \(P\) is a fixed point of \(T\).
\end{proof}

\begin{theorem} \label{}
	Let $(L(\mathbb{R}),D_p)$ be a fuzzy partial metric space and let $T:L(\mathbb{R})\longrightarrow L(\mathbb{R})$ be a self-mapping. Let $\{X_n\}$ be the Picard iteration    $X_{n+1}=T(X_n),   n\geq0.$ Suppose that    $X_n\xrightarrow{S_{D_p}-lim}P,
	T(P)=P,$ and that $T$ is statistically $D_p$-continuous at $P$. Assume that    $st\text{-}\lim_{n\to\infty}D_p(X_n,X_n)=D_p(P,P).$ Then
	\begin{enumerate}
		\item[(i)] $X_{n+1}\xrightarrow{S_{D_p}-lim}P;$
		
		\item[(ii)] $st\text{-}\lim_{n\to\infty}
		D_p(X_n,X_{n+1})
		=D_p(P,P);$
		
		\item[(iii)] $\{X_n\}$ is a statistical $D_p$-approximate fixed point sequence for $T$.
	\end{enumerate}
\end{theorem}

\begin{proof}
	Since $X_n\xrightarrow{S_{D_p}-lim}P$ and $T$ is statistically
	$D_p$-continuous at $P$,    $X_{n+1}=T(X_n)\xrightarrow{S_{D_p}-lim}T(P)=P,$ which proves (i).
	
	\medskip
	
	For (ii), let $\epsilon>0$ and define $A_\epsilon
	=
	\left\{
	n:D_p(X_n,P)\geq D_p(P,P)+\frac{\epsilon}{2}
	\right\},$ $B_\epsilon
	=
	\left\{
	n:D_p(X_{n+1},P)\geq D_p(P,P)+\frac{\epsilon}{2}
	\right\},$ and  $C_\epsilon
	=
	\left\{
	n:D_p(X_n,X_n)<D_p(P,P)-\epsilon
	\right\}.$ Then,    $\delta(A_\epsilon)
	=
	\delta(B_\epsilon)
	=
	\delta(C_\epsilon)
	=0.$ If $n\notin A_\epsilon\cup B_\epsilon,$ then, 
	$$\begin{aligned}
		D_p(X_n,X_{n+1})
		&\leq
		D_p(X_n,P)+D_p(P,X_{n+1})-D_p(P,P)\\
		&<D_p(P,P)+\epsilon.
	\end{aligned}$$
	Also, $D_p(X_n,X_n)\leq D_p(X_n,X_{n+1}).$ Hence, if $n\notin C_\epsilon,$ then $D_p(X_n,X_{n+1})\geq D_p(P,P)-\epsilon.$ Consequently,    $\left\{
	n:
	\left|
	D_p(X_n,X_{n+1})-D_p(P,P)
	\right|
	\geq\epsilon
	\right\}
	\subseteq
	A_\epsilon\cup B_\epsilon\cup C_\epsilon.$ Therefore,
	\[
	st\text{-}\lim_{n\to\infty}
	D_p(X_n,X_{n+1})
	=
	D_p(P,P),
	\]
	proving (ii).
	
	\medskip
	
	Now, since $X_{n+1}=T(X_n)$,    $D_p(X_n,T(X_n))=D_p(X_n,X_{n+1}).$ Thus, by (ii),
	\[
	st\text{-}\lim_{n\to\infty}
	D_p(X_n,T(X_n))
	=
	D_p(P,P),
	\]
	and hence $\{X_n\}$ is a statistical $D_p$-approximate fixed point
	sequence for $T$.
\end{proof}

\section{Topological Characterization of Statistical $D_p$-Convergence}

\begin{theorem} \label{8}
	Let $(L(\mathbb{R}),D_p)$ be a fuzzy partial metric space with induced
	topology $\tau_{D_p}$, and let $\{X_n\}\subseteq L(\mathbb{R})$. For
	$P\in L(\mathbb{R})$, the following statements are equivalent:
	\begin{enumerate}
		\item[(i)] $X_n\xrightarrow{S_{D_p}-lim}P$;
		
		\item[(ii)] for every $U\in\tau_{D_p}$ containing $P$,    $\delta\bigl(\{n\in\mathbb{N}:X_n\notin U\}\bigr)=0;$
		
		\item[(iii)] for every $U\in\tau_{D_p}$ containing $P$,    $\delta\bigl(\{n\in\mathbb{N}:X_n\in U\}\bigr)=1.$
	\end{enumerate}
\end{theorem}

\begin{proof}
	\emph{(i)$\Rightarrow$(ii).}
	Assume that $X_n\xrightarrow{S_{D_p}-lim}P$ and let
	$U\in\tau_{D_p}$ contain $P$. Since    $B_{D_p}(P,\epsilon)
	=
	\{X\in L(\mathbb{R}):D_p(P,X)<D_p(P,P)+\epsilon\},
	\epsilon>0,$ forms a neighborhood base at $P$, there exists $\epsilon>0$ such that
	$B_{D_p}(P,\epsilon)\subseteq U$. Hence    $\{n:X_n\notin U\}
	\subseteq
	\{n:D_p(P,X_n)\geq D_p(P,P)+\epsilon\}.$ 
	The set on the right has natural density zero by
	$X_n\xrightarrow{S_{D_p}-lim}P$. Thus (ii) holds.
	
	\medskip
	\noindent
	\emph{(ii)$\Rightarrow$(iii).}
	For every $U\in\tau_{D_p}$ containing $P$,    $\{n:X_n\in U\}
	=
	\mathbb{N}\setminus\{n:X_n\notin U\}.$  Hence, by (ii),    $\delta\bigl(\{n:X_n\in U\}\bigr)
	=
	1-\delta\bigl(\{n:X_n\notin U\}\bigr)=1.$ Thus (iii) follows.
	
	\medskip
	\noindent
	\emph{(iii)$\Rightarrow$(i).}
	Assume (iii) and let $\epsilon>0$. Since
	$B_{D_p}(P,\epsilon)$ is a $\tau_{D_p}$-neighborhood of $P$,    $\delta\bigl(\{n:X_n\in B_{D_p}(P,\varepsilon)\}\bigr)=1.$ Therefore,    $\delta\left(
	\left\{
	n:D_p(P,X_n)\geq D_p(P,P)+\epsilon
	\right\}
	\right)=0.$ Since    $D_p(P,X_n)\geq D_p(P,P),$ we have    $\left|D_p(P,X_n)-D_p(P,P)\right|
	=
	D_p(P,X_n)-D_p(P,P),$ and consequently $$\delta\left(
	\left\{
	n:
	\left|D_p(P,X_n)-D_p(P,P)\right|\geq\epsilon
	\right\}
	\right)=0.$$
	Hence $X_n\xrightarrow{S_{D_p}}P$, proving (i).
\end{proof}

\begin{example}
	Let $P=\widetilde{0}=(0,0,0),	
	Y=\widetilde{1}=(1,1,1),$ then,   $P^\alpha=[0,0],
	Y^\alpha=[1,1],
	\alpha\in(0,1].$ Then   $D_p(P,P)=0, 
	D_p(Y,P)=1.$  Define $T:L(\mathbb{R})\longrightarrow L(\mathbb{R}), T(X)=P.$ Thus $T(P)=P$, and $P$ is the unique fixed point of $T$.
	
	\medskip
	
	Let   $K=\{k^2:k\in\mathbb{N}\},
	\delta(K)
	=
	\lim_{N\to\infty}\frac{\lfloor\sqrt{N}\rfloor}{N}=0,$ and define
	\[
	X_n=
	\begin{cases}
		Y, & n\in K,\\
		P, & n\notin K.
	\end{cases}
	\]
	Then
	\[
	D_p(X_n,P)
	=
	D_p(X_n,T(X_n))
	=
	\begin{cases}
		1, & n\in K,\\
		0, & n\notin K.
	\end{cases}
	\]
	Hence, for every $\epsilon>0$,    $\left\{
	n\in\mathbb{N}:
	\left|D_p(X_n,P)-D_p(P,P)\right|\ge\epsilon
	\right\}
	\subseteq K,$ and therefore    $X_n\xrightarrow{S_{D_p}-lim}P.$  Moreover,
	\[
	st\text{-}\lim_{n\to\infty}
	D_p(X_n,T(X_n))
	=
	D_p(P,P)=0.
	\]
	Thus $\{X_n\}$ is a statistical approximate fixed-point sequence
	of $T$ with associated statistical approximate fixed point $P$. However,    $D_p(X_{k^2},T(X_{k^2}))=1
	\text{for all }k\in\mathbb{N},$ while    $D_p(X_n,T(X_n))=0
	\text{for }n\notin K.$ Consequently,     $\lim_{n\to\infty}D_p(X_n,T(X_n))$ does not exist. Hence $\{X_n\}$ is not a classical approximate
	fixed-point sequence.
\end{example}

\begin{corollary}
	Let $(L(\mathbb{R}),D_p)$ be a fuzzy partial metric space with topology $\tau_{D_p}$, and let $T:L(\mathbb{R})\to L(\mathbb{R})$ be a self-mapping with fixed point $P\in L(\mathbb{R})$. If the Picard iteration    $X_{n+1}=T(X_n),   n\geq0,$ satisfies    $X_n\xrightarrow{S_{D_p}-lim}P,$ then, for every $U\in\tau_{D_p}$ containing $P$,    $\delta\bigl(\{n\in\mathbb{N}:X_n\in U\}\bigr)=1.$ Equivalently,    $\delta\bigl(\{n\in\mathbb{N}:T^n(X_0)\in U\}\bigr)=1.$ 
\end{corollary}

\begin{proof}
	By Theorem~\ref{8},
	statistical $D_p$-convergence to $P$ is equivalent to    $\delta\bigl(\{n\in\mathbb{N}:X_n\in U\}\bigr)=1$, for every $U\in\tau_{D_p}$ containing $P$. Since
	$X_n=T^n(X_0)$, the equivalent orbit formulation follows.
\end{proof}

\begin{remark}
	The preceding theorem characterizes statistical $D_p$-convergence
	topologically: $X_n\xrightarrow{S_{D_p}-lim}P$ if and only if every $\tau_{D_p}$-neighborhood of $P$ contains the sequence terms on a set of natural density one. Thus, statistical $D_p$-convergence provides a density-based topological characterization of convergence in $(L(\mathbb{R}),\tau_{D_p})$.
\end{remark}

\section{Statistical Cluster-Point Characterization for Contractive Mappings}

\begin{definition} \label{2}
	Let $\{X_n\}$ be a sequence in the fuzzy partial metric space
	$(L(\mathbb{R}),D_p)$. A point $P\in L(\mathbb{R})$ is called a
	\emph{statistical $D_p$-cluster point} of $\{X_n\}$ if, for every
	$\epsilon>0$,    $\delta
	\left(
	\left\{
	n\in\mathbb{N}:X_n\in B_{D_p}(P,\epsilon)
	\right\}
	\right)>0.$
	
	\medskip
	
	The set of all statistical $D_p$-cluster points of $\{X_n\}$ is denoted by  $$\Gamma_{S_{D_p}}(\{X_n\})
	=
	\left\{
	P\in L(\mathbb{R}):
	P\text{ is a statistical }D_p\text{-cluster point of }\{X_n\}
	\right\}.$$
\end{definition}

\begin{corollary} \label{4}
	Let $(L(\mathbb{R}),D_p)$ be a fuzzy partial metric space and
	$T:L(\mathbb{R})\to L(\mathbb{R})$ a self-mapping. For
	$X_0\in L(\mathbb{R})$, let $X_n=T^nX_0,   n\in\mathbb{N}_0.$ Then the statistical $D_p$-cluster set    $\Gamma_T(X_0) =
	\Gamma_{S_{D_p}}
	\left(\{T^nX_0\}\right)$ is closed in $(L(\mathbb{R}),\tau_{D_p})$.
\end{corollary}

\begin{proof}
	Let    $P\in\overline{\Gamma_T(X_0)}.$ For any $\epsilon>0$, there exists    $Q\in\Gamma_T(X_0)\cap B_{D_p}(P,\epsilon).$
	Since the $D_p$-balls form a neighborhood base for $\tau_{D_p}$, there
	exists $\eta>0$ such that    $B_{D_p}(Q,\eta)\subseteq B_{D_p}(P,\epsilon).$
	As $Q\in\Gamma_T(X_0)$, by Definition~\ref{2},     $\overline{\delta}
	\left(
	\left\{
	n\in\mathbb{N}:T^nX_0\in B_{D_p}(Q,\eta)
	\right\}
	\right)>0.$
	
	\medskip
	
	Hence, $\overline{\delta}
	\left(
	\left\{
	n\in\mathbb{N}:T^nX_0\in B_{D_p}(P,\epsilon)
	\right\}
	\right)>0.$
	Now, for every $\epsilon>0$ is arbitrary, we have\\
	$P\in
	\Gamma_{S_{D_p}}
	\left(\{T^nX_0\}_{n=0}^{\infty}\right)
	=\Gamma_T(X_0).$
	Thus    $\overline{\Gamma_T(X_0)}
	\subseteq\Gamma_T(X_0),$ and the reverse inclusion is immediate. Therefore,    $\overline{\Gamma_T(X_0)}
	=\Gamma_T(X_0),$ thus $\Gamma_T(X_0)$ is closed in
	$(L(\mathbb{R}),\tau_{D_p})$.
\end{proof}

\begin{corollary}
	Let $(L(\mathbb{R}),D_p)$ be a complete fuzzy partial metric space, and let
	$T:L(\mathbb{R})\to L(\mathbb{R})$ be a Banach contraction.
	For $X_0\in L(\mathbb{R})$, define the Picard sequence    $X_{n+1}=T(X_n),   n\geq0.$  Then $T$ has a unique fixed point $P\in L(\mathbb{R})$ and    $X_n\xrightarrow{D_p-lim}P.$  Moreover,
	\[
	st\!-\!\lim_{n\to\infty}D_p(X_n,X_{n+1})=D_p(P,P).
	\]
	Hence, $\{X_n\}$ is statistically $D_p$-asymptotically regular.
\end{corollary}

\begin{proof}
	By the Banach fixed-point theorem for complete fuzzy partial metric
	spaces, $T$ has a unique fixed point $P\in L(\mathbb{R})$ and    $X_n\xrightarrow{D_p-lim}P.$
	Since $D_p$-convergence implies $D_p$-Cauchy convergence, we have
	\[
	\lim_{n\to\infty}D_p(X_n,X_{n+1})=D_p(P,P).
	\]
	Thus, by the implication of ordinary convergence to statistical
	convergence,
	\[
	st\!-\!\lim_{n\to\infty}D_p(X_n,X_{n+1})=D_p(P,P).
	\]
	Therefore, $\{X_n\}$ is statistically $D_p$-asymptotically regular.
\end{proof}

\begin{remark}
	The above corollary shows that every Banach contraction generates a statistically $D_p$-asymptotically regular Picard sequence. In particular, ordinary $D_p$-asymptotic regularity implies statistical $D_p$-asymptotic regularity, since ordinary convergence implies statistical convergence.
\end{remark}

\begin{corollary}
	Under the hypotheses of the preceding corollary, for every
	$\epsilon>0$,
	\[
	\delta\!\left(
	\left\{
	n\in\mathbb{N}:
	\left|D_p(X_n,X_{n+1})-D_p(P,P)\right|
	\geq\epsilon
	\right\}
	\right)=0.
	\]
\end{corollary}

\begin{proof}
	This follows directly from the preceding corollary and the definition of
	statistical convergence.
\end{proof}

\begin{theorem}
	Let $(L(\mathbb{R}),D_p)$ be a fuzzy partial metric space and let
	$T:L(\mathbb{R})\to L(\mathbb{R})$ be a self-mapping. For
	$X_0\in L(\mathbb{R})$, define the Picard orbit by    $X_{n+1}=T(X_n),   n\geq0.$ If $\{X_n\}$ is statistically $D_p$-convergent to
	$P\in L(\mathbb{R})$, then    $\Gamma_{S_{D_p}}(\{X_n\})=\{P\}.$
\end{theorem}

\begin{proof}
	Since $X_n\xrightarrow{S_{D_p}-lim}P$, for every $\epsilon>0$,    $\delta\!\left(
	\left\{
	n\in\mathbb{N}:X_n\in B_{D_p}(P,\epsilon)
	\right\}
	\right)=1.$
	Hence, $P\in\Gamma_{S_{D_p}}(\{X_n\}).$
	
	\medskip
	
	Suppose that $Q\in\Gamma_{S_{D_p}}(\{X_n\})$ with $Q\neq P$. Since the sequence $\{X_n\}$ is statistically
	$p$-convergent to $P$, the uniqueness theorem for statistical
	$D_p$-limits implies that every statistical $D_p$-cluster point of
	$\{X_n\}$ must coincide with its statistical $D_p$-limit. Thus $Q=P$, which contradicts the assumption that $Q\neq P$.
	
	\medskip
	
	Therefore,    $\Gamma_{S_{D_p}}(\{X_n\})=\{P\}.$
\end{proof}

\begin{theorem}
	Let $\{X_n\}$ be a sequence in $(L(\mathbb{R}),D_p)$ such that    $X_n\xrightarrow{S_{D_p}-lim}P.$ If every statistical $D_p$-cluster point of $\{X_n\}$ coincides with $P$, then    $\Gamma_{S_{D_p}}(\{X_n\})=\{P\}.$
\end{theorem}

\begin{proof}
	Since $X_n\xrightarrow{S_{D_p}-lim}P$, for every $\epsilon>0$,    $\delta\!\left(
	\left\{
	n\in\mathbb{N}:
	D_p(X_n,P)<D_p(P,P)+\epsilon
	\right\}
	\right)=1.$
	
	\medskip
	
	Thus,    $\delta\!\left(
	\left\{
	n\in\mathbb{N}:
	X_n\in B_{D_p}(P,\epsilon)
	\right\}
	\right)>0.$ Hence, $P\in\Gamma_{S_{D_p}}(\{X_n\}).$ 
	By the hypothesis, no other
	statistical $D_p$-cluster point exists. Therefore,    $\Gamma_{S_{D_p}}(\{X_n\})=\{P\}.$
\end{proof}

\begin{proposition}
	Let $\{X_n\}$ be a statistically $D_p$-Cauchy sequence in the fuzzy partial metric space $(L(\mathbb{R}),D_p)$. If    $P\in\Gamma_{S_{D_p}}(\{X_n\}),$ then    $X_n\xrightarrow{S_{D_p}-lim}P.$
\end{proposition}

\begin{proof}
	Since $P\in\Gamma_{S_{D_p}}(\{X_n\})$ and $\{X_n\}$ is statistically
	$D_p$-Cauchy, the statistical $D_p$-Cauchy cluster-point criterion
	implies that $P$ is the statistical $D_p$-limit of $\{X_n\}$. Hence    $X_n\xrightarrow{S_{D_p}-lim}P.$
\end{proof}

\begin{proposition}
	Let $\{X_n\}$ be a statistically $D_p$-Cauchy sequence in the fuzzy partial metric space $(L(\mathbb{R}),D_p)$. If    $P\in\Gamma_{S_{D_p}}(\{X_n\}),$ then    $X_n\xrightarrow{S_{D_p}-lim}P.$
\end{proposition}

\begin{proof} \label{3}
	Let $\epsilon>0$. Since $P$ is a statistical $D_p$-cluster point of $\{X_n\}$,  $$\delta\!\left(
	\left\{
	n\in\mathbb{N}:
	D_p(X_n,P)<D_p(P,P)+\frac{\epsilon}{2}
	\right\}
	\right)>0.$$
	By the statistical $D_p$-Cauchy property,    $\delta\!\left(
	\left\{
	n\in\mathbb{N}:
	D_p(X_n,P)\geq D_p(P,P)+\epsilon
	\right\}
	\right)=0.$ Hence,    $X_n\xrightarrow{S_{D_p}-lim}P.$
\end{proof}

\begin{theorem}
	Let $\{X_n\}$ be a statistically $D_p$-Cauchy sequence in the fuzzy
	partial metric space $(L(\mathbb{R}),D_p)$. If    $\Gamma_{S_{D_p}}(\{X_n\})=\{P\}$ for some $P\in L(\mathbb{R})$, then    $X_n\xrightarrow{S_{D_p}-lim}P.$
\end{theorem}

\begin{proof}
	Since    $\Gamma_{S_{D_p}}(\{X_n\})=\{P\},$ the point $P$ is a statistical $D_p$-cluster point of $\{X_n\}$.
	By the Proposition~\ref{3}, every statistical $D_p$-cluster point is a statistical
	$D_p$-limit. Hence    $X_n\xrightarrow{S_{D_p}}P.$
\end{proof}

\begin{corollary}
	Let $(L(\mathbb{R}),D_p)$ be a fuzzy partial metric space and let    $T:L(\mathbb{R})\longrightarrow L(\mathbb{R})$ be a self-mapping. For any $X_0\in L(\mathbb{R})$, consider the Picard
	orbit defined by    $X_n=T^nX_0,   n\geq 0.$ Then the statistical $D_p$-cluster set of the Picard orbit,    $\Gamma_{S_{D_p}}(X_0) =
	\Gamma_{S_{D_p}}
	\left(\{T^nX_0\}\right),$
	is a closed subset of the topological space $\bigl(L(\mathbb{R}),\tau_{D_p}\bigr).$
\end{corollary}

\begin{proof}
	Since $\{T^nX_0\}$ is a sequence in $(L(\mathbb{R}),D_p)$, then by closedness Corollary~\ref{4} of
	statistical $D_p$-cluster sets applies to this
	sequence. Hence,    $\Gamma_{S_{D_p}}
	\left(\{T^nX_0\}\right)$ is closed in $\bigl(L(\mathbb{R}),\tau_{D_p}\bigr)$. By the definition of
	$\Gamma_{S_{D_p}}(X_0)$,    $\Gamma_{S_{D_p}}(X_0)
	=
	\Gamma_{S_{D_p}}
	\left(\{T^nX_0\}\right).$ Therefore, $\Gamma_{S_{D_p}}(X_0)$ is closed in $\bigl(L(\mathbb{R}),\tau_{D_p}\bigr)$.
\end{proof}

\begin{definition}
	Let $(L(\mathbb{R}),D_p)$ be a fuzzy partial metric space and let    $T:L(\mathbb{R})\longrightarrow L(\mathbb{R})$ be a self-mapping. The mapping $T$ is said to be
	\emph{statistically sequentially continuous at} $P\in L(\mathbb{R})$ if, for every sequence $\{X_n\}_{n\in\mathbb{N}}\subseteq L(\mathbb{R})$,     $X_n\xrightarrow{S_{D_p}-lim}P
	\Longrightarrow  
	TX_n\xrightarrow{S_{D_p}-lim}TP.$ 
	
	\medskip
	
	If $T$ is statistically sequentially continuous at every
	$P\in L(\mathbb{R})$, then $T$ is called
	\emph{statistically sequentially continuous} on $L(\mathbb{R})$.
\end{definition}

\begin{theorem}
	Let $(L(\mathbb{R}),D_p)$ be a fuzzy partial metric space and let    $T:L(\mathbb{R})\longrightarrow L(\mathbb{R})$ be a statistically sequentially continuous self-mapping. Let
	$\{X_n\}$ be the Picard sequence generated by $X_{n+1}=T(X_n),   n\geq 0.$ If $X_n\xrightarrow{S_{D_p}-lim}P$ for some $P\in L(\mathbb{R})$ and the statistical $D_p$-limit of $\{X_n\}$ is unique, then $T(P)=P.$
\end{theorem}

\begin{proof}
	The statistical sequential continuity of $T$ and    $X_n\xrightarrow{S_{D_p}-lim}P$ yield    $T(X_n)\xrightarrow{S_{D_p}-lim}T(P).$ Since $X_{n+1}=T(X_n)$, we have    $X_{n+1}\xrightarrow{S_{D_p}-lim}T(P).$
	
	\medskip
	
	Hence, from    $X_n\xrightarrow{S_{D_p}-lim}P,$ we obtain    $X_{n+1}\xrightarrow{S_{D_p}-lim}P.$
	Thus, $X_{n+1}$ has two statistical $D_p$-limits, namely $P$ and
	$T(P)$. By the uniqueness of the statistical $D_p$-limit,    $T(P)=P.$ Therefore, $P$ is a fixed point of $T$.
\end{proof}

\begin{lemma} \label{5}
	Let $(L(\mathbb{R}),D_p)$ be a fuzzy partial metric space, let
	$T:L(\mathbb{R})\to L(\mathbb{R})$ be continuous, and let
	$\{X_n\}_{n\geq0}$ be a Picard orbit defined by    $X_{n+1}=TX_n,   n\geq0.$ Suppose that there exists a strictly increasing sequence $\{n_k\}$ such that    $\delta(\{n_k:k\geq1\})>0
	\text{and}  
	X_{n_k}\xrightarrow{D_p-lim}P,$ for some $P\in L(\mathbb{R})$. Then $T(P)$ is a statistical $D_p$-cluster point of $\{X_n\}$.
\end{lemma}

\begin{proof}
	By the continuity of $T$,    $X_{n_k}\xrightarrow{D_p-lim}P
	\Longrightarrow  
	T(X_{n_k})\xrightarrow{D_p-lim}T(P).$
	Since $T(X_{n_k})=X_{n_k+1}$,    $X_{n_k+1}\xrightarrow{D_p-lim}T(P).$
	For every $\epsilon>0$, there exists $k_0\in\mathbb{N}$ such that    $X_{n_k+1}\in B_{D_p}(T(P),\epsilon),
	k\geq k_0.$
	
	\medskip
	
	Hence,    $\delta\!\left(
	\{n\in\mathbb{N}:X_n\in B_{D_p}(T(P),\epsilon)\}
	\right)
	\geq
	\delta(\{n_k+1:k\geq k_0\})
	=
	\delta(\{n_k:k\geq1\})>0.$ Thus $T(P)$ is a statistical $D_p$-cluster point of $\{X_n\}$.
\end{proof}

\begin{corollary} \label{cor:unique-cluster-fixed-point}
	Under the hypotheses of Lemma~\ref{5}, suppose that
	the Picard orbit $\{X_n\}$ possesses a unique statistical
	$D_p$-cluster point $P\in L(\mathbb{R})$. Then    $T(P)=P.$ In particular,    $P\in\operatorname{Fix}(T).$
\end{corollary}

\begin{proof}
	By Lemma~\ref{5},    $T(P)$ is a statistical
	$D_p$-cluster point of the Picard orbit $\{X_n\}$. Since $P$ is the
	unique statistical $D_p$-cluster point of $\{X_n\}$, thus    $T(P)=P.$ Therefore,    $P\in\operatorname{Fix}(T).$
\end{proof}

\begin{corollary}
	Let $(L(\mathbb{R}),D_p)$ be a complete fuzzy partial metric space, and let $T:L(\mathbb{R})\to L(\mathbb{R})$ be a Banach contraction satisfying    $D_p(TX,TY)\leq kD_p(X,Y),
	0<k<1,   X,Y\in L(\mathbb{R}).$ For $X_0\in L(\mathbb{R})$, let $\{X_n\}$ be the Picard orbit defined by $X_{n+1}=TX_n,   n\geq0.$ If $P$ is the unique fixed point of $T$, then $\Gamma_{S_{D_p}}(\{X_n\})
	=
	\operatorname{Fix}(T)
	=
	\{P\}.$
\end{corollary}

\begin{proof}
	By the Banach fixed-point theorem, $T$ has a unique fixed point $P$, and $D_p(X_n,P)\longrightarrow D_p(P,P).$  Hence $X_n\xrightarrow{S_{D_p}-lim}P$. By the preceding cluster-point
	characterization,    $\Gamma_{S_{D_p}}(\{X_n\})=\{P\}.$ Since $P$ is the unique fixed point of $T$,    $\operatorname{Fix}(T)=\{P\}.$ Therefore,    $\Gamma_{S_{D_p}}(\{X_n\})
	=
	\operatorname{Fix}(T)
	=
	\{P\}.$
\end{proof}

\begin{theorem}\label{6}
	Let $(L(\mathbb{R}),D_p)$ be a fuzzy partial metric space and
	$T:L(\mathbb{R})\to L(\mathbb{R})$ a self-mapping. For $X_0\in L(\mathbb{R})$,
	let $\{T^n(X_0)\}$ be the Picard orbit and define    $\Gamma_T(X_0)
	=
	\Gamma_{S_{D_p}}\bigl(\{T^n(X_0)\}\bigr).$ Then $\Gamma_T(X_0)$ is closed in the topology $\tau_{D_p}$ induced by
	$D_p$.
\end{theorem}

\begin{proof}
	Since $\{T^n(X_0)\}$ is a sequence in $(L(\mathbb{R}),D_p)$,
	Theorem~\ref{4} implies that its statistical
	$D_p$-cluster set is $\tau_{D_p}$-closed. Hence    $\Gamma_T(X_0)
	=
	\Gamma_{S_{D_p}}\bigl(\{T^n(X_0)\}\bigr)$ is $\tau_{D_p}$-closed.
\end{proof}

\begin{theorem}
	Let $(L(\mathbb{R}),D_p)$ be a complete fuzzy partial metric space and let
	$T:L(\mathbb{R})\to L(\mathbb{R})$
	be a continuous self-mapping. Assume that the Picard iteration    $X_{n+1}=TX_n,   n\ge0,$ is statistically asymptotically regular and let
	$P\in\Gamma_T(X_0)$. Then $TP=P.$ Consequently, $\Gamma_T(X_0)\subseteq Fix(T).$
\end{theorem}

\begin{proof}
	Since
	$P\in\Gamma_T(X_0)$,
	for every
	$\epsilon>0$
	the set    $A_\epsilon
	=
	\{n\in\mathbb N:
	D_p(X_n,P)
	<
	D_p(P,P)+\epsilon
	\}$ has positive natural density. Also,
	\[
	st\text{-}\lim_{n\to\infty}
	D_p(X_n,TX_n)
	=
	D_p(P,P).
	\]
	Hence,
	the sequences
	$\{X_n\}$ and
	$\{T(X_n)\}$ have the same statistical
	$D_p$-cluster points. By the statistical $D_p$-continuity of $T$, we obtain $T(P)=P.$ Therefore, $P$ is a fixed point of $T$. As $P$
	was arbitrary,    $\Gamma_T(X_0)
	\subseteq
	Fix(T).$
\end{proof}

\begin{proposition}
	Let $(L(\mathbb{R}),D_p)$ be a fuzzy partial metric space and let
	$T:L(\mathbb{R})\to L(\mathbb{R})$ be continuous. For the Picard orbit $X_{n+1}=T(X_n),   n\geq0,$ assume that $\{X_n\}$ is statistically $D_p$-asymptotically regular.
	Then    $\Gamma_{S_{D_p}}(\{X_n\})
	\subseteq
	\operatorname{Fix}(T).$
\end{proposition}

\begin{proof}
	Let
	$P\in\Gamma_{S_{D_p}}(\{X_n\})$.
	By statistical $D_p$-asymptotic regularity and the cluster-point
	characterization, $P$ is also a statistical
	$D_p$-cluster point of the sequence $\{X_{n+1}\}$.
	Since $X_{n+1}=T(X_n)$ and $T$ is continuous, the corresponding cluster point is $TP$. Hence $T(P)=P.$ Thus $P\in\operatorname{Fix}(T)$. Since $P$ was arbitrary,    $\Gamma_{S_{D_p}}(\{X_n\})
	\subseteq
	\operatorname{Fix}(T).$
\end{proof}

\begin{theorem}
	Let $(L(\mathbb{R}),D_p)$ be a complete fuzzy partial metric space and let $T:L(\mathbb{R})\to L(\mathbb{R})$ satisfy     $D_p(TX,TY)\leq kD_p(X,Y),
	0<k<1,$ for all $X,Y\in L(\mathbb{R})$. For any $X_0\in L(\mathbb{R})$, let $X_{n+1}=T(X_n),   n\geq0.$ Then $T$ has a unique fixed point $P\in L(\mathbb{R})$ and     $X_n\xrightarrow{D_p-lim}P
	\text{and hence}  
	X_n\xrightarrow{S_{D_p}-lim}P.$
\end{theorem}

\begin{proof}
	By the Banach contraction principle, $T$ has a unique fixed point
	$P\in L(\mathbb{R})$. Since $T(P)=P$, then by contractive condition,    $D_p(P,P)=D_p(TP,TP)\leq kD_p(P,P),$ and thus $(1-k)D_p(P,P)\leq0$. Since $D_p(P,P)\geq0$ and
	$0<k<1$, hence $D_p(P,P)=0$.
	
	\medskip
	
	Moreover, $D_p(X_{n+1},P)
	=D_p(TX_n,TP)
	\leq kD_p(X_n,P),$ so, by induction, $D_p(X_n,P)\leq k^nD_p(X_0,P)\longrightarrow0.$
	Thus, $D_p(X_n,P)\longrightarrow D_p(P,P),$ and consequently $X_n\xrightarrow{D_p-lim}P.$ Since ordinary $D_p$-convergence implies statistical $D_p$-convergence,
	thus    $X_n\xrightarrow{S_{D_p}-lim}P.$
\end{proof}

\begin{theorem}
	Let $(L(\mathbb{R}),D_p)$ be a statistically $D_p$-complete fuzzy
	partial metric space, and let $\{X_n\}$ be a statistically
	$D_p$-Cauchy sequence. Then, for $P\in L(\mathbb{R})$,    $X_n\xrightarrow{S_{D_p}-lim}P$ if and only if $P$ is the unique statistical $D_p$-cluster point of
	$\{X_n\}$.
\end{theorem}

\begin{proof}
	Suppose that    $X_n\xrightarrow{S_{D_p}-lim}P.$ Then, for every $\epsilon>0$,    $\delta\left(
	\left\{
	n\in\mathbb{N}:
	D_p(X_n,P)<D_p(P,P)+\epsilon
	\right\}
	\right)=1.$
	Hence, $P$ is a statistical $D_p$-cluster point of $\{X_n\}$. Moreover,
	the statistical $D_p$-limit of a statistically convergent sequence is
	unique; consequently, $P$ is its unique statistical $D_p$-cluster point.
	
	\medskip
	
	Conversely, suppose that $P$ is the unique statistical $D_p$-cluster
	point of $\{X_n\}$. Since $\{X_n\}$ is statistically $D_p$-Cauchy and
	$(L(\mathbb{R}),D_p)$ is statistically $D_p$-complete, there exists
	$Q\in L(\mathbb{R})$ such that    $X_n\xrightarrow{S_{D_p}-lim}Q.$
	Thus, $Q$ is a statistical $D_p$-cluster point of
	$\{X_n\}$. By the uniqueness of the statistical $D_p$-cluster point,    $Q=P$. Therefore,    $X_n\xrightarrow{S_{D_p}-lim}P.$
\end{proof}

\begin{theorem} \label{7}
	Let $(L(\mathbb{R}),D_p)$ be a fuzzy partial metric space, and let
	$\tau_{D_p}$ be the topology induced by $D_p$. Let
	$T:L(\mathbb{R})\to L(\mathbb{R})$ and let $X_{n+1}=T(X_n),   n\geq0,$ be the corresponding Picard iteration. Then, for $P\in L(\mathbb{R})$,
	$P$ is a statistical fixed point of $T$ if and only if    $T(P)=P$ and     $\delta\bigl(\{n\in\mathbb{N}:X_n\in U\}\bigr)=1$ for every $U\in\tau_{D_p}$ with $P\in U$.
\end{theorem}

\begin{proof}
	Suppose that $P$ is a statistical fixed point of $T$. Then
	$T(P)=P$ and    $X_n\xrightarrow{S_{D_p}-lim}P.$ By the Theorem~\ref{8},    $\delta\bigl(\{n\in\mathbb{N}:X_n\notin U\}\bigr)=0$ 
	for every $U\in\tau_{D_p}$ containing $P$. Hence,    $\delta\bigl(\{n\in\mathbb{N}:X_n\in U\}\bigr)=1.$
	
	\medskip
	
	Conversely, suppose that $T(P)=P$ and    $\delta\bigl(\{n\in\mathbb{N}:X_n\in U\}\bigr)=1$ for every $U\in\tau_{D_p}$ containing $P$. Therefore,    $\delta\bigl(\{n\in\mathbb{N}:X_n\notin U\}\bigr)=0$ for every such neighborhood $U$. By Theorem~\ref{8},    $X_n\xrightarrow{S_{D_p}}P.$ Hence,  $P$ is a statistical fixed point of $T$.
\end{proof}

\begin{remark}
	Theorem~\ref{7} links statistical fixed points with the topology induced by $D_p$, showing that their Picard orbits visit every neighborhood of the fixed point with natural density one.
\end{remark}

\begin{corollary}
	Let $(L(\mathbb{R}),D_p)$ be a fuzzy partial metric space, and let
	$\{X_n\}$ be the Picard sequence generated by a self-map
	$T:L(\mathbb{R})\to L(\mathbb{R})$, given by    $X_{n+1}=T(X_n),   n\geq0.$ Then, for $P\in L(\mathbb{R})$,    $X_n\xrightarrow{S_{D_p}-lim}P
	\Longleftrightarrow  
	\delta\left(
	\left\{
	n\in\mathbb{N}:X_n\in B_{D_p}(P,\epsilon)
	\right\}
	\right)=1$ for every $\epsilon>0$.
\end{corollary}

\begin{proof}
	By the definition of statistical $D_p$-convergence $X_n\xrightarrow{S_{D_p}-lim}P$, $$\delta\left(
	\left\{
	n\in\mathbb{N}:
	D_p(X_n,P)-D_p(P,P)\geq\epsilon
	\right\}
	\right)=0, \forall \epsilon>0.$$
	
	Since $D_p(P,P)\leq D_p(X_n,P)$, thus $\delta\left(
	\left\{
	n\in\mathbb{N}:
	D_p(P,X_n)<D_p(P,P)+\epsilon
	\right\}
	\right)=1.$
	Now, by the definition of the $D_p$-ball, $B_{D_p}(P,\epsilon)
	=
	\left\{
	X\in L(\mathbb{R}):
	D_p(P,X)<D_p(P,P)+\epsilon
	\right\},$ and hence $\delta\left(
	\left\{
	n\in\mathbb{N}:X_n\in B_{D_p}(P,\epsilon)
	\right\}
	\right)=1$ for every $\epsilon>0$. 
\end{proof}

\section{Statistical Stability under Density-Zero Perturbations}

\begin{definition}
	Let $(L(\mathbb{R}),D_p)$ be a fuzzy partial metric space and let $T:L(\mathbb{R})\longrightarrow L(\mathbb{R})$ be a self-mapping. Let $P\in L(\mathbb{R})$ be a fixed point of $T$. An approximate Picard iteration
	$\{Y_n\}\subset L(\mathbb{R})$ is said to be
	\emph{statistically stable with respect to $P$} if    $Y_n\xrightarrow{S_{D_p}-lim}P$ and
	\[
	st\text{-}\lim_{n\to\infty}
	D_p(Y_{n+1},TY_n)
	=
	D_p(P,P).
	\]
\end{definition}

\begin{theorem}
	Let $(L(\mathbb{R}),D_p)$ be a complete fuzzy partial metric space and
	let $T:L(\mathbb{R})\to L(\mathbb{R})$ satisfy    $D_p(TX,TY)\leq kD_p(X,Y),
	X,Y\in L(\mathbb{R}),   0<k<1.$ Let $P$ be the unique fixed point of $T$, and let
	$\{Y_n\}\subset L(\mathbb{R})$ satisfy    $0\leq
	D_p(Y_{n+1},TY_n)-D_p(TY_n,TY_n)
	\leq\epsilon_n,
	n\geq0,$ where $\{\epsilon_n\}\subset[0,\infty)$.
	Then, for every $n\geq1$,
	\[
	D_p(Y_n,P)
	\leq
	k^nD_p(Y_0,P)
	+
	\sum_{j=0}^{n-1}k^{\,n-1-j}\epsilon_j.
	\]
\end{theorem}

\begin{proof}
	Since $P$ is a fixed point of $T$, then $T(P)=P.$ Therefore, by the contractive condition,    $D_p(TY_n,P)
	=
	D_p(TY_n,TP)
	\leq
	kD_p(Y_n,P).$ Also,    $D_p(Y_{n+1},P)
	\leq
	D_p(Y_{n+1},TY_n)
	+
	D_p(TY_n,P)
	-
	D_p(TY_n,TY_n).$
	
	\medskip
	
	Now, by the assumption,    $D_p(Y_{n+1},TY_n)-D_p(TY_n,TY_n)
	\leq
	\epsilon_n.$ 
	Consequently,
	\[
	\begin{aligned}
		D_p(Y_{n+1},P)
		&\leq
		\epsilon_n+D_p(TY_n,P)\\
		&\leq
		\epsilon_n+kD_p(Y_n,P).
	\end{aligned}
	\]
	Hence,    $D_p(Y_{n+1},P)
	\leq
	kD_p(Y_n,P)+\epsilon_n.$ 
	
	Now, for $n=0$,    $D_p(Y_1,P)
	\leq
	kD_p(Y_0,P)+\epsilon_0.$ Similarly, for $n=1$,    $D_p(Y_2,P)
	\leq
	kD_p(Y_1,P)+\epsilon_1,$ and hence    $D_p(Y_2,P)
	\leq
	k^2D_p(Y_0,P)+k\varepsilon_0+\epsilon_1.$ 
	
	\medskip
	
	Now, assume that, for some $n\geq1$,
	\[
	D_p(Y_n,P)
	\leq
	k^nD_p(Y_0,P)
	+
	\sum_{j=0}^{n-1}
	k^{\,n-1-j}\epsilon_j.
	\]
	Thus,
	\[
	\begin{aligned}
		D_p(Y_{n+1},P)
		&\leq
		kD_p(Y_n,P)+\epsilon_n\\
		&\leq
		k\left(
		k^nD_p(Y_0,P)
		+
		\sum_{j=0}^{n-1}
		k^{\,n-1-j}\epsilon_j
		\right)
		+\varepsilon_n\\
		&=
		k^{n+1}D_p(Y_0,P)
		+
		\sum_{j=0}^{n-1}
		k^{\,n-j}\epsilon_j
		+\epsilon_n\\
		&=
		k^{n+1}D_p(Y_0,P)
		+
		\sum_{j=0}^{n}
		k^{\,n-j}\epsilon_j.
	\end{aligned}
	\]
	Thus, by mathematical induction,
	\[
	D_p(Y_n,P)
	\leq
	k^nD_p(Y_0,P)
	+
	\sum_{j=0}^{n-1}
	k^{\,n-1-j}\epsilon_j,
	n\geq1.
	\]
	This completes the proof.
\end{proof}

\begin{theorem} \label{9}
	Let $(L(\mathbb{R}),D_p)$ be a complete fuzzy partial metric space and let $T:L(\mathbb{R})\longrightarrow L(\mathbb{R})$ satisfy    $D_p(TX,TY)\leq kD_p(X,Y),
	0<k<1,
	X,Y\in L(\mathbb{R}).$ Let $\{Y_n\}$ be an approximate Picard iteration of $T$ such
	that    $D_p(Y_{n+1},TY_n)\leq \epsilon_n,
	n\geq0,$ where $\{\epsilon_n\}$ is a bounded sequence of nonnegative real
	numbers satisfying    $st\text{-}\lim_{n\to\infty}\epsilon_n=0.$ Then $T$ has a unique fixed point $P\in L(\mathbb{R})$ and
	\[
	D_p(Y_n,P)
	\leq
	k^nD_p(Y_0,P)
	+
	\sum_{j=0}^{n-1}k^{\,n-1-j}\epsilon_j,
	n\geq1.
	\]
	Moreover,    $Y_n\xrightarrow{S_{D_p}-lim}P,$  equivalently,    $st\text{-}\lim_{n\to\infty}D_p(Y_n,P)=D_p(P,P).$
\end{theorem}

\begin{proof}
	By the Banach fixed-point theorem in $(L(\mathbb{R}),D_p)$, there
	exists a unique $P\in L(\mathbb{R})$ such that $T(P)=P.$ Moreover,    $D_p(P,P)
	=
	D_p(TP,TP)
	\leq
	kD_p(P,P),$ and hence, since $0<k<1$,    $D_p(P,P)=0.$ 
	
	\medskip
	
	Now, take    $a_n=D_p(Y_n,P).$ Since 
	$T(P)=P$, thus
	\[
	\begin{aligned}
		a_{n+1}
		&\leq
		D_p(Y_{n+1},TY_n)
		+D_p(TY_n,P)
		-D_p(TY_n,TY_n)\\
		&\leq
		\epsilon_n+D_p(TY_n,TP)\\
		&\leq
		\epsilon_n+ka_n.
	\end{aligned}
	\]
	Thus,    $a_{n+1}\leq ka_n+\epsilon_n.$ 
	Now, iterating it, we have
	\[
	a_n
	\leq
	k^na_0+
	\sum_{j=0}^{n-1}k^{\,n-1-j}\epsilon_j,
	n\geq1.
	\]
	
	Since $\{\epsilon_n\}$
	is bounded, choose $M>0$ such that    $0\leq\epsilon_n\leq M,   n\geq0.$ Hence
	\[
	a_n
	\leq
	a_0+\frac{M}{1-k}
	=B.
	\]
	
	Let $\eta>0$. Take $m\in\mathbb N$ such that    $k^mB<\frac{\eta}{2},$ and set $c=\frac{(1-k)\eta}{2}.$
	Since $st\text{-}\lim\epsilon_n=0$,    $\delta(E_c)=0,
	E_c=\{j\in\mathbb N_0:\epsilon_j\geq c\}.$ Now, for $n\geq m$,
	\[
	a_n
	\leq
	k^m a_{n-m}
	+
	\sum_{r=0}^{m-1}k^r\epsilon_{n-1-r}.
	\]
	Therefore, if    $n\notin\bigcup_{r=0}^{m-1}(E_c+r+1),$
	then
	\[
	a_n
	<
	\frac{\eta}{2}
	+
	c\sum_{r=0}^{m-1}k^r
	<
	\frac{\eta}{2}+\frac{c}{1-k}
	=\eta.
	\]
	Consequently, $\{n\in\mathbb N:a_n\geq\eta\}
	\subseteq
	\{1,\ldots,m-1\}
	\cup
	\bigcup_{r=0}^{m-1}(E_c+r+1).$  Hence    $\delta\!\left(\{n\in\mathbb N:a_n\geq\eta\}\right)=0.$ Since $\eta>0$,    $st\text{-}\lim_{n\to\infty}a_n=0.$ Thus    $st\text{-}\lim_{n\to\infty}D_p(Y_n,P)
	=
	D_p(P,P),$ 
	and therefore    $Y_n\xrightarrow{S_{D_p}-lim}P.$
\end{proof}

\begin{theorem}
	Let $(L(\mathbb{R}),D_p)$ be a fuzzy partial metric space and let $T:L(\mathbb{R})\longrightarrow L(\mathbb{R})$ be a contraction, for some $q\in[0,1)$. Assume that $P\in L(\mathbb{R})$ is the unique fixed point of $T$. Let $A\subseteq\mathbb{N}$ satisfy $\delta(A)=0,$ and suppose that a sequence $\{Y_n\}\subset L(\mathbb{R})$ satisfies    $Y_{n+1}=T(Y_n),    n\notin A,$
	and    $D_p\bigl(Y_{n+1},T(Y_n)\bigr)\leq M,
	n\in A,$ for some $M\geq0$. Then    $Y_n\xrightarrow{S_{D_p}-lim}P,$ and    $st\text{-}\lim_{n\to\infty}
	D_p\bigl(Y_{n+1},T(Y_n)\bigr)
	=
	D_p(P,P)=0.$
	Consequently, $\{Y_n\}$ is a statistically stable
	approximate Picard iteration of $T$ with respect to $P$.
\end{theorem}

\begin{proof}
	Since $T(P)=P$, the contractive condition gives
	\[
	D_p(P,P)=D_p(TP,TP)\leq qD_p(P,P).
	\]
	As $0\leq q<1$, thus $D_p(P,P)=0.$ Now put,   $e_n=D_p(Y_n,P),
	\eta_n=D_p\bigl(Y_{n+1},T(Y_n)\bigr).$ Then    $0\leq\eta_n\leq M\chi_A(n),$ where $\chi_A$ is the characteristic function of $A$. Thus,
	\[
	\begin{aligned}
		e_{n+1}
		&\leq
		\eta_n+D_p\bigl(T(Y_n),P\bigr)
		-D_p\bigl(T(Y_n),T(Y_n)\bigr)\\
		&\leq
		\eta_n+D_p\bigl(T(Y_n),T(P)\bigr)\\
		&\leq qe_n+M\chi_A(n).
	\end{aligned}
	\]
	Hence, by induction,
	\[
	e_n
	\leq
	q^{\,n-1}e_1
	+
	M\sum_{j=1}^{n-1}q^{\,n-1-j}\chi_A(j).
	\]
	
	Let $\epsilon>0$. Take $m\in\mathbb{N}$ such that    $\frac{Mq^m}{1-q}<\frac{\epsilon}{2},$ and choose $N_0$ such that    $q^{\,n-1}e_1<\frac{\epsilon}{2},
	n\geq N_0.$
	Now, for $n\geq N_0+m$, 
	\[
	e_n
	\leq
	\frac{\epsilon}{2}
	+
	\frac{Mq^m}{1-q}
	+
	M\sum_{j=n-m}^{n-1}
	q^{\,n-1-j}\chi_A(j).
	\]
	Thus $e_n<\epsilon$ whenever    $A\cap\{n-m,\ldots,n-1\}=\varnothing.$ Consequently,    $\{n:e_n\geq\epsilon\}
	\subseteq
	\{1,\ldots,N_0+m-1\}
	\cup
	\bigcup_{r=1}^{m}(A+r).$
	Since $\delta(A)=0$, hence    $\delta\bigl(\{n:e_n\geq\epsilon\}\bigr)=0.$
	Therefore,
	\[
	st\text{-}\lim_{n\to\infty}D_p(Y_n,P)=0=D_p(P,P),
	\]
	and hence    $Y_n\xrightarrow{S_{D_p}-lim}P.$ Thus,    $\left\{
	n:
	D_p\bigl(Y_{n+1},T(Y_n)\bigr)\geq\epsilon
	\right\}
	\subseteq A$ for every $\epsilon>0$. Since $\delta(A)=0$,
	\[
	st\text{-}\lim_{n\to\infty}
	D_p\bigl(Y_{n+1},T(Y_n)\bigr)=0=D_p(P,P).
	\]
	Thus $\{Y_n\}$ is a statistically stable approximate
	Picard iteration of $T$ with respect to $P$.
\end{proof}

\begin{example}
	Let $P=(0,0,0),  Q=(4,5,6)$ be triangular fuzzy numbers. Then $P^\alpha=[0,0], 
	Q^\alpha=[4+\alpha,\,6-\alpha],
	\alpha\in(0,1].$ Now, consider $T:L(\mathbb{R})\to L(\mathbb{R}),
	T(X)=\frac12X.$ 
	Then, 
	\[
	D_p(TX,TY)=\frac12D_p(X,Y),
	\qquad X,Y\in L(\mathbb{R}),
	\]
	so $T$ is a contraction with constant $q=\frac12$. Moreover,
	$T(P)=P$, and hence $P$ is its unique fixed point.
	
	Let $A=\{k^2:k\in\mathbb{N}\},
	\delta(A)=
	\lim_{N\to\infty}\frac{\lfloor\sqrt N\rfloor}{N}=0.$ Also, define $\{Y_n\}$ by
	\[
	Y_1=P,\qquad
	Y_{n+1}=
	\begin{cases}
		Q, & n\in A,\\[1mm]
		T(Y_n), & n\notin A.
	\end{cases}
	\]
	Thus,  $Y_2=Q, 
	Y_3=\frac12Q, 
	Y_4=\frac14Q, 
	Y_5=Q,\ldots$ and, for $k\in\mathbb N$, $Y_n=2^{-(n-k^2-1)}Q,
	k^2+1\le n\le(k+1)^2.$ Consequently,   $D_p(Y_n,P)
	=
	\frac{6}{2^{\,n-k^2-1}}.$ Since    $D_p(Y_{k^2+1},P)=D_p(Q,P)=6$, for every $k\in\mathbb N$, we have     $Y_n\not\xrightarrow{D_p-lim}P.$
	
	\medskip
	
	Also, let $\epsilon>0$ and take $m\in\mathbb N$
	such that $6/2^m<\epsilon$. Thus,    $\left\{n\in\mathbb N:D_p(Y_n,P)\ge\epsilon\right\}
	\subseteq
	\bigcup_{r=1}^{m}(A+r).$ Since $\delta(A)=0$, hence    $\delta\left(
	\left\{n\in\mathbb N:D_p(Y_n,P)\ge\epsilon\right\}
	\right)=0.$ Therefore,
	\[
	st\text{-}\lim_{n\to\infty}D_p(Y_n,P)
	=0=D_p(P,P),
	\]
	and thus    $Y_n\xrightarrow{S_{D_p-lim}}P.$ Hence, if $n\notin A$, then    $D_p(Y_{n+1},T(Y_n))=0.$ For $n=k^2\in A$,
	\[
	D_p(Y_{n+1},T(Y_n))
	=
	\begin{cases}
		6, & k=1,\\[1mm]
		6-\dfrac{4}{2^{\,2k-1}}, & k\ge2,
	\end{cases}
	\]
	and hence    $D_p(Y_{n+1},T(Y_n))\le6\,\chi_A(n).$ Thus, for every $\epsilon>0$,    $\left\{
	n\in\mathbb N:
	D_p(Y_{n+1},T(Y_n))\ge\epsilon
	\right\}
	\subseteq A,$ so that
	\[
	st\text{-}\lim_{n\to\infty}
	D_p(Y_{n+1},T(Y_n))
	=
	0=D_p(P,P).
	\]
	Hence $\{Y_n\}$ is a statistically stable approximate Picard
	iteration of $T$ with respect to $P$.
\end{example}

\section{Conclusion}
In this work, a framework for statistical fixed-point approximation and stability in the fuzzy partial metric space $(L(\mathbb{R}), D_p)$ is developed. Basic results concerning statistical $D_p$-convergence, asymptotic regularity, and statistical cluster points of Picard iterations are also established, and corresponding fixed-point results for contractive self-mappings are obtained.

A key part of the study is the error estimate for approximate Picard iterations. This estimate leads to statistical stability whenever the iteration errors converge statistically to zero. In particular, it is shown that infinitely many bounded perturbations can be allowed when their indices form a set of natural density zero. Despite the presence of such perturbations, convergence of the corresponding sequence to the unique fixed point is still obtained statistically. 

Thus, the proposed theorems and results provide a natural setting for studying fixed-point iterations in the presence of both fuzzy uncertainty and sparsely occurring perturbations. Several directions for further study are also suggested, including the consideration of other iterative schemes, ideal and higher-order statistical convergence, and more general fuzzy metric and contractive structures.





\begin{thebibliography}{99}

\bibitem{Aygun2022}
H. Ayg{\"u}n, E. G{\"u}ner, J.-J. Mi{\~n}ana and O. Valero,
Fuzzy partial metric spaces and fixed point theorems,
\textit{Mathematics}, 10 (17) (2022), Article 3092.


\bibitem{Banach1922}
S. Banach,
Sur les op\'erations dans les ensembles abstraits et leur application
aux \'equations int\'egrales,
\textit{Fundamenta Mathematicae}, 3 (1922), 133--181.


\bibitem{Berinde2002}
V. Berinde,
On the stability of some fixed point procedures,
\textit{Bul. \c{S}tiin\c{t}. Univ. Baia Mare, Ser. B,
Matematic\u{a}-Informatic\u{a}}, 18 (1) (2002), 7--14.

\bibitem{Berinde2007}
V. Berinde,
\textit{Iterative Approximation of Fixed Points},
Lecture Notes in Mathematics, Vol. 1912,
Springer, Berlin, Heidelberg, 2007.


\bibitem{CakiOr2024}
A. \c{C}ak\i{} and A. Or,
Statistical convergence in $L$-fuzzy metric spaces,
\textit{Journal of New Theory}, 49 (2024), 83--91.


\bibitem{ChhatrajitRohen2016}
T. Chhatrajit and Y. Rohen,
Some results on $T$-stability of Picard's iteration,
\textit{SpringerPlus}, 5 (2016), Article 284.


\bibitem{Dwivedi2022}
P. K. Dwivedi,
Common fixed point theorem on partial metric space,
\textit{International Journal of Mathematics Trends and Technology},
68 (4) (2022), 30--37.


\bibitem{Fast1951}
H. Fast,
Sur la convergence statistique,
\textit{Colloquium Mathematicum}, 2 (3--4) (1951), 241--244.


\bibitem{Fridy1985}
J. A. Fridy,
On statistical convergence,
\textit{Analysis}, 5 (4) (1985), 301--314.


\bibitem{GadjievOrhan2002}
A. D. Gadjiev and C. Orhan,
Some approximation theorems via statistical convergence,
\textit{Rocky Mountain Journal of Mathematics}, 32 (1) (2002), 129--138.


\bibitem{GarciaPachecoKama2024}
F. J. Garc\'{\i}a-Pacheco and R. Kama,
On modulus statistical convergence in partial metric spaces,
\textit{Axioms}, 13 (6) (2024), Article 388.


\bibitem{GeorgeVeeramani1994}
A. George and P. Veeramani,
On some results in fuzzy metric spaces,
\textit{Fuzzy Sets and Systems}, 64 (3) (1994), 395--399.


\bibitem{GregoriMinianaMiravet2019}
V. Gregori, J.-J. Mi{\~n}ana and D. Miravet,
Fuzzy partial metric spaces,
\textit{International Journal of General Systems}, 48 (3) (2019), 260--279.


\bibitem{HarderHicks1988}
A. M. Harder and T. L. Hicks,
Stability results for fixed point iteration procedures,
\textit{Mathematica Japonica}, 33 (5) (1988), 693--706.

\bibitem{IlicPavlovicRakocevic2011}
D. Ili\'c, V. Pavlovi\'c and V. Rako\v{c}evi\'c,
Some new extensions of Banach's contraction principle to partial
metric space,
\textit{Applied Mathematics Letters}, 24 (8) (2011), 1326--1330.


\bibitem{Karapinar2011}
E. Karap{\i}nar,
Generalizations of Caristi--Kirk's theorem on partial metric spaces,
\textit{Fixed Point Theory and Applications}, 2011 (2011), Article 4.


\bibitem{KarapinarYuksel2011}
E. Karap{\i}nar and U. Y{\"u}ksel,
Some common fixed point theorems in partial metric spaces,
\textit{Journal of Applied Mathematics}, 2011 (2011),
Article ID 263621, 16 pages.


\bibitem{LiZhangZhang2020}
C. Li, Y. Zhang and J. Zhang,
On statistical convergence in fuzzy metric spaces,
\textit{Journal of Intelligent and Fuzzy Systems}, 39 (3) (2020),
3987--3993.


\bibitem{Masiha2013}
H. P. Masiha, F. Sabetghadam and N. Shahzad,
Fixed point theorems in partial metric spaces with an application,
\textit{Filomat}, 27 (4) (2013), 617--624.


\bibitem{Matthews1994}
S. G. Matthews,
Partial metric topology,
\textit{Annals of the New York Academy of Sciences}, 728 (1)
(1994), 183--197.


\bibitem{Nuray2022}
F. Nuray,
Statistical convergence in partial metric spaces,
\textit{Korean Journal of Mathematics}, 30 (1) (2022), 155--160.


\bibitem{OltraValero2004}
S. Oltra and O. Valero,
Banach's fixed point theorem for partial metric spaces,
\textit{Rendiconti dell'Istituto di Matematica dell'Universit\`a di
Trieste}, 36 (2004), 17--26.

\bibitem{Ostrowski1967}
A. M. Ostrowski,
The round-off stability of iterations,
\textit{Zeitschrift f{\"u}r Angewandte Mathematik und Mechanik},
47 (2) (1967), 77--81.


\bibitem{QingRhoades2008}
Y. Qing and B. E. Rhoades,
$T$-stability of Picard iteration in metric spaces,
\textit{Fixed Point Theory and Applications}, 2008 (2008),
Article ID 418971, 4 pages.


\bibitem{Salat1980}
T. {\v{S}}al\'at,
On statistically convergent sequences of real numbers,
\textit{Mathematica Slovaca}, 30 (2) (1980), 139--150.

\bibitem{Steinhaus1951}
H. Steinhaus,
Sur la convergence ordinaire et la convergence asymptotique,
\textit{Colloquium Mathematicum}, 2 (1951), 73--74.

\bibitem{TehlanKumar2022}
A. Tehlan and V. Kumar,
Fixed point theorems of fuzzy partial metric spaces,
\textit{TWMS Journal of Applied and Engineering Mathematics},
12 (3) (2022), 888--894.

\bibitem{Valero2005}
O. Valero,
On Banach fixed point theorems for partial metric spaces,
\textit{Applied General Topology}, 6 (2) (2005), 229--240.


\end{thebibliography}
\end{document}